\documentclass[11pt]{article}

\usepackage{makeidx}
\usepackage{amssymb}
\usepackage{amsfonts}
\usepackage{amsmath}
\usepackage{euscript}
\usepackage{theorem}
\usepackage{enumerate}

\theorembodyfont{\rm}      
\theoremstyle{change}      
\begingroup
\newtheorem{thm}{Theorem\hskip 5mm}[section]

\newtheorem{prop}[thm]{Proposition\hskip 5mm}

\newtheorem{cor}[thm]{Corollary\hskip 5mm}

\newtheorem{lem}[thm]{Lemma\hskip 5mm}

\endgroup
\begingroup

\endgroup
\begingroup

\endgroup
\newcommand{\qed}{{\unskip\nobreak\hfil\penalty 50\hskip 2em\hbox{}
 \nobreak\hfil$\square$\parfillskip=0pt\finalhyphendemerits=0\par}}

\newenvironment{interlist}%
{\begin{list}{}{%
\setlength{\labelsep}{0pt}\setlength{\leftmargin}{-2.5em}%
\setlength{\labelwidth}{0pt}%
\setlength{\listparindent}{0pt}}}%
{\end{list}}
{\begin{interlist}%
     \vspace{-.5em}\item}
{\end{interlist}}

\def\X{{X^\pm}}
\def\E{{ E^\pm}}
\def\a{{\alpha}}
\def\b{{\beta}}
\def\k{{\kappa}}
\def\t{{\tau}}
\def\e{{\epsilon}}
\def\l{{\lambda}}
\def\S{{\widehat {S}}}
\def\f{{\widetilde {f}}}
\def\Lt{{\widetilde {L}}}
\def\LH{{\widehat{L}}}

\def\Sp{\mathrm {Sp}}
\def\GL{\mathrm {GL}}
\def\ker{{\rm ker}\,}
\def\ind{{\rm ind}}
\def\Q{{\mathbf Q}}

\begin{document}

\begin{center}
On the 2-modular reduction of the Steinberg representation of the
symplectic group
\end{center}
\begin{center}
Fernando Szechtman
\end{center}
\begin{center}
{\small\textit{Department of Pure Mathematics, University of
Waterloo, Ontario, Canada, N2L 3G1\\ e-mail:
fszechtm@herod.uwaterloo.ca }}
\end{center}

\begin{abstract}

We show that in characteristic 2, the Steinberg representation of
the symplectic group $\Sp_{2n}(q)$, $q$ a power of an odd prime
$p$, has two irreducible constituents lying just above the socle
that are isomorphic to the two Weil modules of degree $(q^n-1)/2$.

\end{abstract}

\section{Introduction}

The Steinberg representation was constructed by R. Steinberg
\cite{St} in 1957 for all finite Chevalley groups $G=G(F_q)$, $q$
a power of a prime $p$. He did so by identifying a particular
right ideal in the group algebra $FG$ over an an arbitrary field
$F$. He then showed this ideal to be minimal precisely when the
characteristic of $F$ does not divide the index of the normalizer
of a $p$-Sylow subgroup of $G$. In particular, the Steinberg
representation is irreducible in characteristic 0.

The general problem of finding the composition factors of the
Steinberg module has recently been addressed by Gow \cite{G}, who
has a conjecture in the case of the general linear group. He also
determined all irreducible constituents in a few examples. One
composition factor that is certainly known is the socle, which is
irreducible and generated as a $G$-module by the fixed points of
$U$.

Our main result in this context is the following: in
characteristic 2, the Steinberg representation of the symplectic
group $\Sp_{2n}(q)$, $q$ a power of an odd prime $p$, has two
irreducible constituents lying just above the socle that are
isomorphic to the two Weil modules of degree $(q^n-1)/2$.

We start the paper by constructing a family of irreducible
representations for a $p$-Sylow subgroup $U$ of $\Sp_{2n}(q)$ over
a field $F$ of characteristic $l\neq p$ containing a primitive
$p$-th root of unity. This is done in section 3 under no
restrictions on $p$. We also produce families of irreducible
$F$-modules for the normalizer $B$ of $U$, and for a maximal
parabolic subgroup $P$ of $\Sp_{2n}(q)$. It is then shown in
section 4 that, for $p$ odd, the above irreducible modules are
essentially those found in the Weil representation. All
information needed in regards to the Weil representation is
contained in section 4.

In section 5 we show that, provided $l$ divides $q+1$, the
Steinberg module restricted to $P$ contains copies of the two Weil
modules of degree $(q^n-1)/2$ restricted to $P$, up to
multiplication by an explicit linear character $P\to\{\pm
1\}\subset F^*$.

Finally, in section 6 we prove that, if $l=2$, the Steinberg
module modulo the socle contains copies of the two Weil modules
of degree $(q^n-1)/2$. In order to extend the results of section 5
to those of section 6 we prove a criterion ensuring that only two
identities in the symplectic group algebra $F\Sp_{2n}(q)$ need to
be verified. The actual verification of these identities is
difficult, and occupies the rest of the paper.

The most basic concepts and definitions, along with our choice of
notation can be found in section 2.

I am very grateful to R. Gow for suggesting this problem.









\section{Basic Notions}

\noindent {\bf The Symplectic Group.} Let $F_q$ be a finite field
of characteristic $p$. Let $V$ be a vector space of dimension $2n$
over $F_q$ endowed with a non-degenerate alternating form
$\langle ~,~\rangle$. The symplectic group of rank $2n$ over
$F_q$, denoted simply by $\Sp$, is the group of all $g\in \GL(V)$
satisfying
$$\langle gv,gw\rangle=\langle v,w\rangle,\quad v,w\in V.$$

\noindent {\bf Witt decomposition of $V$ and associated subgroups
of $\Sp$.} We fix throughout a symplectic basis of $V$. This is a
basis $u_1,...,u_n,v_1,...,v_n$ satisfying
$$
\langle u_i,v_j\rangle=\delta_{ij},\quad \langle
u_i,u_j\rangle=0,\quad \langle v_i,v_j\rangle=0.
$$
Write $M=\mathrm{span}\,\{u_1,...,u_n\}$ and
$N=\mathrm{span}\,\{v_1,...,v_n\}$. Let $\Sp^M$ be the pointwise
stabilizer of $M$ in $\Sp$. The matrix of any $s\in \Sp^M$ has
the appearance
\begin{equation}
\label{ese} \left(\begin{array}{c|c} 1 & S \\\hline 0 & 1
\end{array}\right),
\end{equation}
where $S$ is a symmetric $n\times n$ matrix with coefficients in
$F_q$. Let $\Sp_{M,N}$ be the subgroup of $\Sp$ preserving both
$M$ and $N$. The matrix of any $a\in\Sp_{M,N}$ has the form
\begin{equation}
\label{a} \left(\begin{array}{c|c} A & 0 \\\hline 0 & {}^t A^{-1}
\end{array}\right),
\end{equation} where $A\in\GL_n(q)$. The group $\Sp_{M,N}$
normalizes $\Sp^M$.
The semidirect product group $\Sp_M=\Sp^M\rtimes \Sp_{M,N}$ is the
global stabilizer of $M$ in $\Sp$. It is a maximal subgroup of
$\Sp$.

\noindent {\bf A $p$-Sylow subgroup of $\Sp$.} Let $T$ be the
group of all $g$ in $\Sp_{M,N}$ whose matrix is of the form
(\ref{a}), where $A$ is upper triangular and has 1's on the main
diagonal. Then $U=T\rtimes \Sp^M$ is a $p$-Sylow subgroup of
$\Sp$.

Let $H$ be the group of all $g\in\Sp_{M,N}$ having all $u_i$ and
$v_j$ as eigenvalues. The group $U$ is normalized by $H$. The
semidirect product $B=U\rtimes H$ is the normalizer of $U$ in
$\Sp$.

\noindent {\bf The Weyl group.} For $1\leq i< n$, denote by $w_i$
the element of $\Sp$ that preserves the given symplectic basis,
and is defined by
$$
w_i=(u_i,u_{i+1})(v_{i},v_{i+1}).
$$
We also define, for $1\leq j\leq n$, the elements $c_j$ of $\Sp$
by: $c_j(u_j)=v_j$, $c_j(v_j)=-u_j$, while all other basis vectors
remain fixed.

Let $W_0=\langle w_1,...,w_{n-1}\rangle$ and $W_1=\langle
c_1,...,c_{n}\rangle$. The group $W_1$ is normalized by $W_0$. Let
$W_2=W_1\rtimes W_0$. The group $H$ is normalized by $W_2$, with
$W_2\cap H=W_1\cap H=\langle c_1^2,...,c_n^2\rangle$. The Weyl
group of $\Sp$ is
$$
W=W_2H/H\cong W_2/W_2\cap H\cong W_1\rtimes W_0/W_1\cap H\cong
(W_1/W_1\cap H)\rtimes W_0.
$$
It is the split extension of $C_2^n\cong \langle
c_1,...,c_n\rangle /\langle c_1^2,...,c_n^2\rangle$ by $S_n\cong
\langle w_1,...,w_{n-1}\rangle$. Here the symmetric group $S_n$
acts naturally on $C_2^n$ by permuting the $n$ copies of $C_2$.

Set ${\mathcal {F}}=\{ w_1,...,w_{n-1},c_n\}$. Given $w\in W_2$,
write $\ell(w)=\ell(wH)$ for the length of the shortest word in
the letters $w_1H$,...,$w_{n-1}H$,$c_n H$ which is equal to $wH$
in the Weyl group $W$. The largest value of $\ell$ is attained by
$w_0=c_1\cdots c_n$, the element of $\Sp$ defined by
$$
w_0(u_i)=v_i,\quad w_0(v_i)=-u_i,\quad 1\leq i\leq n.
$$
Observe that $\ell(w_0)=n^2$.

For $w\in W_2$, set $U_w^+= w^{-1}Uw \cap U$, the group of all
$u\in U$ whose conjugates $w uw^{-1}$ remain in $U$. Let
$U_w^{-}$ be the the group of all $u\in U$ whose conjugates $w
uw^{-1}$ belong to $w_0Uw_0^{-1}$. One has $\Sp=U^+_wU^-_w$ with
$U^+_w\cap U^{-}_w=1$. Furthermore, $|U_w^-|=q^{\ell(w)}$,
$|U_w^+|=q^{n^2-\ell(w)}$ and $|U|=q^{n^2}$.

\noindent {\bf Root subgroups.} Set ${\mathcal
{I}}=\{(i,j)\,\vert\, 1\leq i<j\leq 2n,\; i\leq n\}$. Let
$(i,j)\in {\mathcal {I}}$ and take $\alpha\in F_q$. If $j\leq n$
define $x_{i,j}(\alpha)$ to be the element of $T$ whose matrix
has the form (\ref{a}), where $A_{i,j}=\alpha$, $A_{kk}=1$ for
all $k$, and all other entries of $A$ are equal to zero. If $j>n$
define $x_{ij}(\alpha)$ to be the element of $\Sp^M$ whose matrix
has the form (\ref{ese}), where $S_{i,j-n}=S_{j-n,i}=\alpha$, and
all other entries of $S$ are equal to zero. For a fixed $(i,j)\in
{\mathcal {I}}$ write $X_{i,j}$ for the group of all
$x_{i,j}(\alpha)$ with $\alpha\in F_q$. This is the root subgroup
corresponding to $(i,j)$. The fundamental root subgroups are
$X_{(i,i+1)}$,  for $1\leq i<n$, and $X_{(n,2n)}$. Observe that
$X_{(i,i+1)}=U_{w_i}^-$ and $X_{(n,2n)}=U_{c_n}^-$.

\noindent {\bf Symplectic transvections.} Given $u$ in $V$ and
$\alpha\in F_q$, let $\rho_{u,\alpha}\in\Sp$ be the symplectic
transformation defined by
$$
\rho_{u,\alpha}(v)=v+\alpha\langle u,v\rangle u,\quad v\in V.
$$
\noindent {\bf The center of $\Sp$.} A symplectic transformation
commuting with all symplectic transvections is necessarily a
scalar operator. Thus the center $Z(\Sp)$ of $\Sp$ is equal to
$\{1,-1\}$ if $p$ is odd, and is trivial otherwise.

\noindent {\bf Conjugation.} For group elements $g$ and $h$ we
write ${}^h g=hgh^{-1}$ and $g^h=h^{-1}gh$.

\section{A representation of $\Sp_M$}

We fix throughout a field $F$ of characteristic $l\neq p$
possessing a non-trivial $p$-th root of unity. We also fix a
non-trivial linear character $\lambda:F_q^+\to F^*$.

For $v\in V$, let $\chi_v:\Sp\to F^*$ be the function defined by
\begin{equation}
\label{formu} \chi_v(g)=\lambda(\langle gv,v\rangle),\quad g\in
\Sp.
\end{equation}
We shall use the same notation for $\chi_v$ and its restriction
to various subgroups. Context will dictate what subgroup is being
used at a given time.

Given $v$ in $N$, let $S_v$ be the group of all $g$ in $\Sp_M$
satisfying $gv\equiv v\mod M$. From the matrix representation
(\ref{ese}) of $\Sp^M$ we see that $\Sp^M\subset S_v$. In fact,
$S_v=\Sp^M\rtimes (\Sp_{M,N})_v$, where $(\Sp_{M,N})_v$ is the
pointwise stabilizer of $v$ in $\Sp_{M,N}$.

\begin{lem}
\label{homo} Let $v\in N$. Then $\chi_v:S_v\to F^*$ is a  linear
character of $S_v$.
\end{lem}

\noindent{\it Proof.} Let $g$ and $h$ be elements of $S_v$. Then
$$
\begin{aligned}
\langle ghv,v\rangle & = \langle g(hv-v+v),v\rangle \\
& = \langle g(hv-v),v\rangle +  \langle gv,v\rangle \\
& = \langle hv-v,g^{-1}v\rangle +  \langle gv,v\rangle \\
& = \langle hv-v,g^{-1}v-v+v\rangle +  \langle gv,v\rangle.
\end{aligned}
$$
As $h,g^{-1}\in S_v$, we have $hv-v\in M$ and $g^{-1}v-v\in M$,
whence $\langle hv-v,g^{-1}v-v\rangle=0$. It follows that
$$
\langle ghv,v\rangle=\langle hv-v,v\rangle+\langle gv,v\rangle=
\langle hv,v\rangle +  \langle gv,v\rangle.
$$
\qed

For $v\in N$, let $\S_v=S_v\times Z(\Sp)$, the group of all $g$ in
$\Sp_M$ satisfying $gv\equiv \pm v\mod M$. We shall denote by
$\chi_v^\pm$ the group homomorphism $\S_v\to F^*$ that agrees with
$\chi_v$ on $S_v$ and maps $-1\in Z(\Sp)$ to $\pm 1\in F^*$.

Since $\Sp^M$ is contained in $S_v$ for all $v\in N$, we may view
$\chi_v$ as a linear character of $\Sp^M$. A fundamental property
is the following.
\begin{lem}
\label{fp} Let $v,w\in N$. Then
$$
\chi_v=\chi_w\text{ on }\Sp^M\Leftrightarrow v=\pm w.
$$
\end{lem}

\noindent{\it Proof.} One implication follows directly from
(\ref{formu}). Suppose conversely that
\begin{equation}
\label{uno} \lambda(\langle gv,v\rangle)=\lambda(\langle
gw,w\rangle)
\end{equation}
for all $g\in \Sp^M$. Fix any $u\in V$. We apply (\ref{uno}) to
$g=\rho_{u,\alpha}$ obtaining
$$
\lambda(\a(\langle u,v\rangle^2-\langle u,w\rangle^2))=1
$$
for all $\alpha\in F_q$. As $\lambda$ is non-trivial, it follows
that
\begin{equation}
\label{dos} \langle u,v\rangle^2=\langle u,w\rangle^2
\end{equation} for all $u\in V$. We infer that the linear maps
$\langle-,v\rangle$ and $\langle-,w\rangle$ have the same kernel.
Choose $u\in V$ satisfying $\langle u,v\rangle=1$. Then
(\ref{dos}) yields $\langle u,v\rangle=\pm 1$. We conclude that
$\langle-,v\rangle=\langle-,w\rangle$ or
$\langle-,v\rangle=\langle-,-w\rangle$, that is $v=w$ or $v=-w$,
as claimed.\qed

In what follows it will sometimes be convenient to think of the
matrix of an arbitrary $g\in\Sp_{M,N}$ as having the form
\begin{equation}
\label{b} \left(\begin{array}{c|c} {}^t A^{-1} & 0 \\\hline 0 & A
\end{array}\right),
\end{equation} where $A\in\GL_n(q)$. In this case $g\in T$
whenever $A$ is lower triangular and has 1's on the main diagonal.
In particular, $gv_n=v_n$ for all $g\in T$. It follows that $U$
is contained in $S_{\a v_n}$ for all $\a\in F_q$, a fact to be
used repeatedly below.

\begin{lem}
\label{puo} Let $\a\in F_q^*$. Then the linear character
$\chi_{\a v_n}$ of $U$ is trivial on $U_{c_n}^+$ and non-trivial
on $U_{c_n}^-$.
\end{lem}

\noindent{\it Proof.} We have $U_{c_n}^+=(\Sp^M)_0\rtimes T$,
where $(\Sp^M)_0$ is the subgroup of $\Sp^M$ of all $g$ whose
matrix has the form (\ref{ese}) with $S_{nn}=0$. If $g\in T$ then
$gv_n=v_n$. Therefore (\ref{formu}) yields $T\subset\ker \chi_{\a
v_n}$. If $g\in (\Sp^M)_0$ then $gv_n-v_n\in
\mathrm{span}\,\{u_1,...,u_{n-1}\}$. As $v_n$ is orthogonal to
itself and $u_1,...u_{n-1}$, (\ref{formu}) yields
$(\Sp^M)_0\subset \ker \chi_{\a v_n}$, as well. This proves the
first assertion. As for the second, note that $U_{c_n}^-$ is the
group of all $x_{n,2n}(\b)=\rho_{u_n,\b}$, $\b\in F_q$. Now
$$
\chi_{\a v_n}(\rho_{u_n,\b})=\lambda(\langle \a v_n+\b\langle
u_n,\a v_n\rangle u_n,\a v_n\rangle)= \lambda(\b\a^2)
$$
for all $\b\in F_q$. As $\lambda$ is non-trivial, $\chi_{\a v_n}$
is non-trivial on $U_{c_n}^-$.\qed

Since the fundamental root subgroups $X_{i,i+1}$ are contained in
$U_{c_n}^+$ we obtain the following result.

\begin{cor}
\label{rooto} The linear character $\chi_{\a v_n}$, $\a\in F_q^*$,
is non-trivial in exactly one fundamental root subgroup, namely
$X_{(n,2n)}=U^-_{c_n}$.
\end{cor}


\begin{lem}
\label{g0} Let $v\in V$ and $g_0\in\Sp$. Then $$
\chi_v(g^{g_0})=\chi_{g_0v}(g),\quad g\in\Sp.
$$
\end{lem}

\noindent{\it Proof.} Since $g_0$ preserves $\langle~,~\rangle$,
we have
$$
\chi_v(g^{g_0})=\l(\langle g_0^{-1}gg_0v,v\rangle)=\lambda(\langle
gg_0v,g_0v\rangle)=\chi_{g_0v}(g),\quad g\in\Sp.
$$
\qed

Let $X=X_\lambda$ be a vector space over $F$ possessing a basis
$(E_v)_{v\in N}$ indexed by $N$.
\begin{prop}
The function $R:\Sp_M\to\GL(X)$ given by
\begin{equation}
\label{rep} R(sa)E_v=\lambda(\langle
sav,av\rangle)E_{av}=\chi_{av}(s)E_{av},\quad v\in N,s\in \Sp^M,
a\in \Sp_{M,N}
\end{equation}
defines an $F$-representation of $\Sp_M$.
\end{prop}

\noindent{\it Proof.} Suppose $s_1,s_2\in\Sp^M$, $a_1,a_2\in
\Sp_{M,N}$, and $v\in N$. Then
$$
R(s_1a_1s_2a_2)E_v=R(s_1({}^{a_1}s_2)a_1a_2)E_v=
\chi_{a_1a_2v}(s_1({}^{a_1}s_2))E_{a_1a_2v}.$$ Since
$\chi_{a_1a_2v}$ is a linear character of $\Sp^M$
$$
\chi_{a_1a_2v}(s_1({}^{a_1}s_2))=\chi_{a_1a_2v}(s_1)\chi_{a_1a_2v}({}^{a_1}s_2).
$$
By Lemma \ref{g0} we have
$$\chi_{a_1a_2v}({}^{a_1}s_2)=\chi_{a_1a_2v}(s_2^{a_1^{-1}})=\chi_{a_2v}(s_2).$$
On the other hand
$$
R(s_1a_1)R(s_2a_2)E_v=R(s_1a_1)\chi_{a_2v}(s_2)E_{a_2v}=\chi_{a_1a_2v}(s_1)\chi_{a_2v}(s_2)
E_{a_1a_2v}.
$$
This completes the proof.\qed

We determine the irreducible constituents of this representation,
depending on whether $p$ is even or odd. Suppose first that $p$
is odd. For $v\in N\setminus\{0\}$, set $E^+=E_v+E_{-v}$ and
$E^-_v=E_v-E_{-v}$. We also write $X^+=X^+_\lambda$ for the span
of the $E_v^+$, and let $X^-=X^-_\lambda$ be the span of the
$E_v^-$. When $p=2$ we set $E^+_v=E^-_v=E_v$ for all $v\in
N\setminus\{0\}$, and let $X^+=X^+_\lambda=X^-=X^-_\lambda$ be
the span of the $E_v$, $v\neq 0$.

Formula (\ref{rep}) ensures that $\X$ is an $F\Sp_M$-module. Its
dimension is $(q^n-1)/2$ for $p$ odd and $q^n-1$ for $p=2$. We
write $X_0$ for the $\Sp_M$-fixed points of $X$, that is
$X_0=FE_0$.


For $1\leq i\leq n$ and $\a\in F_q^*$, consider the following
$F$-subspace of $\X$:

$$
\X_{\a,i}=\text{span} \{\E_{\a v_i+\a_{i+1}v_{i+1}+\cdots+\a_n
v_n}\,\vert\, \a_i \in F_q\}.
$$

\begin{thm}
\label{irru} $\X_{\a,i}$ is an absolutely irreducible
$FU$-submodule of $\X$ of dimension $q^{n-i}$. Furthermore, if
$w$ is the element of $W_0$ defined by
$$
w=(u_i,u_n,u_{n-1},...,u_{i+1})(v_i,v_n,v_{n-1},...,v_{i+1})
$$
then
$$ X^+_{\a,i}\cong\mathrm{ind}_{U^+_w}^U F E^+_{\a v_i}\cong \mathrm{ind}_{U^+_w}^U F E^-_{\a v_i}\cong X^-_{\a,i},$$
where $F \E_{\a v_i}$ affords the linear character $\chi_{\a
v_i}$ of $U^+_w$.
\end{thm}

\noindent{\it Proof.} Let $\E_{\a v_i+\a_{i+1}v_{i+1}+\cdots+\a_n
v_n}$ be a typical basis vector of $\X_{a,i}$, and let $g\in U$.
Then $g=sa$, where the matrix of $a$ has the form (\ref{b}) with
$A$ lower triangular with 1's on the main diagonal, and $s\in
\Sp^M$.

We see from (\ref{rep})  how $s$ and $a$ act on $\X$. Indeed,
$\E_{\a v_i+\a_{i+1}v_{i+1}+\cdots+\a_nv_n}$ is an eigenvector for
$s$ acting on $\X$, and the above matrix description of $a$
ensures that $a\E_{\a v_i+\a_{i+1}v_{i+1}+\cdots+\a_nv_n}$ belongs
to $\X_{\a,i}$. It follows that $g$ preserves $\X_{\a,i}$, which
is then an $FU$-module. It is clear that $\dim \X_{\a,i}=q^{n-i}$.

Since $\mathrm{char}\, F=l$ is coprime to $|U|=q^{n^2}$, the
$FU$-module $\X_{\a,i}$ is completely reducible. We proceed to
show that the commuting ring of $\X_{\a,i}$ is comprised entirely
by scalar operators. This implies that $\X_{\a,i}$ is absolutely
irreducible.

Let $C$ be an $F$-endomorphism of $\X_{\a,i}$ commuting with the
action of $U$. The group $U_w^-$ consists of all $g\in\Sp_{M,N}$
whose matrix has the form (\ref{b}), where $A$ is a lower
triangular matrix with 1's on the diagonal and all columns
different from column $i$ have zero entries below the diagonal.
Let us write $0,0,...,0,1,\b_{i+1},...,\b_n$ for the entries in
the $i$-th column of such matrix $A$. Here 1 is in the $i$-th
position. Then
$$
g\E_{\a v_i}=\E_{\a v_i+\a\b_{i+1}v_{i+1}+\cdots+\a\b_nv_n}.
$$ It follows
that $U_w^-$ acts transitively on the basis vectors $\E_{\a
v_i+\a_{i+1}v_{i+1}+\cdots+\a_nv_n}$ of $\X_{\a,i}$. Since the
actions of $C$ and $U_w^-$ on these vectors commute, we infer that
$C$ acts diagonally on them.

Given two different basis vectors $\E_{\a
v_i+\a_{i+1}v_{i+1}+\cdots+\a_nv_n}$ and $\E_{\a
v_i+\b_{i+1}v_{i+1}+\cdots+\b_nv_n}$, we have
$$
\a v_i+\a_{i+1}v_{i+1}+\cdots+\a_nv_n\neq \pm (\a
v_i+\b_{i+1}v_{i+1}+\cdots+\b_nv_n).
$$
It follows from Lemma \ref{fp} that for some $s\in \Sp^M$ one has
$$
\chi_{\a v_i+\a_{i+1}v_{i+1}+\cdots+\a_nv_n}(s)\neq \chi_{\a
v_i+\b_{i+1}v_{i+1}+\cdots+\b_nv_n}(s).
$$
Since the diagonal actions of $\Sp^M$ and $C$ on the basis vectors
$\E_{\a v_i+\a_{i+1}v_{i+1}+\cdots+\a_nv_n}$ commute, we deduce
that $C$ must be a scalar operator, as claimed.

From (\ref{rep}) we see that $\X_{\a,i}$ affords a monomial
representation of $FU$. The stabilizer of $F\E_{\a v_i}$ is the
group of all $g\in U$ such that $g=sa$, $s\in\Sp^M$, $a\in T$, and
the $i$-th column of $A$ in the matrix representation (\ref{b}) of
$a$ is equal to the $i$-th column of the identity matrix.
Therefore, the stabilizer of $F\E_{\a v_i}$ is equal to $U^+_w$.
As $U^-_w$ is a left transversal for $U^+_w$ in $U$ acting
transitively on the one dimensional subspaces $F\E_{\a
v_i+\a_{i+1}v_{i+1}+\cdots+\a_nv_n}$, it follows that
$\X_{\a,i}\cong \ind_{U^+_w}^U F\E_{\a v_i}$. Furthermore, $U^+_w$
acts on $F\E_{\a v_i}$ by means of the linear character $\chi_{\a
v_i}$. This completes the proof of the theorem.\qed

\begin{thm}
\label{isu} For a fixed $i$, if $\a,\b\in F_q^*$ satisfy
$\{\a,-\a\}\neq \{\b,-\b\}$ then the $FU$-modules $\X_{\a,i}$ and
$\X_{\b,i}$ are not isomorphic.
\end{thm}

\noindent {\it Proof.} The linear character $\chi_{\a v_i}$ of
$\Sp^M$ enters the restriction of $\X_{\a,i}$ to $\Sp^M$. The
linear characters entering the restriction of $\X_{\b,i}$ to
$\Sp^M$ are of the form $\chi_{\b v_i+\b_{i+1}
v_{i+1}+\cdots+\b_nv_n}$. If $\{\a,-\a\}\neq \{\b,-\b\}$ then none
of them is $\X_{\a,i}$ by Lemma \ref{fp}.\qed

Note that for $h\in H$ we have
\begin{equation}
\label{h} h\X_{\a,i}=\X_{h_i\a,i},
\end{equation}
where $hv_i=h_iv_i$. We shall henceforth denote by ${\cal{ T}}_q$
a transversal for the action of $\{\pm 1\}$ on $F_q^*$. Thus
$|{\cal{ T}}_q|=(q-1)/2$ if $p$ is odd, and $|{\cal{ T}}_q|=q-1$
if $p=2$. Observe that $\X_{\a, i}=\X_{-\alpha,i}$ for all $\a\in
{\cal{ T}}_q$.

\begin{thm}
\label{irb} For a fixed $i$, the direct sum
$\X_i=\oplus_{\a\in{\cal{T}}_q} \X_{\a,i}$ is an absolutely
irreducible $FB$-module.
\end{thm}

\noindent {\it Proof.} Let $Z$ be an $FB$-submodule of $\X_i$.
Then $Z$ is an $FU$-submodule of $\X_i$. It follows from Theorems
\ref{irru} and \ref{isu} that $Z$ must be the direct sum of some
$\X_{\a,i}$. Hence $Z$ is all of $\X_i$ by (\ref{h}). As this
argument works for $F$ or any extension thereof, $\X_{\a,i}$ is
absolutely irreducible.\qed

\begin{thm}
\label{isb} The $\X_i$ are not isomorphic $FB$-modules.
\end{thm}

\noindent {\it Proof.}  Observe that the dimension of $\X_i$ is
equal to $(q-1) q^{n-i}/2$ if $p$ is odd and $(q-1)q^{n-i}$ if
$p=2$. The result thus follows.\qed

The following relation holds for all $w\in W_0$:
\begin{equation}
\label{wes} w E_{v}=E_{w(v)},\quad v\in N.
\end{equation}

\begin{thm}
\label{irrspm} $\X$ is an absolutely irreducible $\Sp_M$-module.
Moreover,
$$
\X=\mathrm{ind}_{\S_{v_n}}^{\Sp_M} F\E_{v_n},
$$
where $\S_{v_n}$ acts on the one-dimensional subspace $F\E_{
v_n}$ by means of $\chi_{v_n}^\pm$.
\end{thm}

\noindent {\it Proof.} Let $Z$ be an $F\Sp_M$-submodule of $\X$.
Then $Z$ is an $FB$-submodule of $\X$. It follows from Theorems
\ref{irb} and \ref{isb} that $Z$ must be the direct sum of some
$\X_i$. Hence at least one $\E_{v_i}$ belongs to $Z$. Therefore
all $\E_{v_j}$ belong to $Z$, by applying (\ref{wes}) with all
$w_k$. As the $\E_{v_j}$ generate $\X_j$ as an $FB$-submodule,
and the sum of all $\X_j$ equals $\X$, it follows that $Z=\X$.
Again, as this argument works for $F$ or any extension thereof,
$\X$ is absolutely irreducible.

Now $\X$ affords a monomial representation of $\Sp_M$; the
stabilizer of $F\E_{v_n}$ is equal to $\S_{v_n}$; the linear
character of $\S_{v_n}$ afforded by $F\E_{v_n}$ is equal to
$\chi_{v_n}^\pm$; the index of $S_{v_n}^+$ in $\Sp_M$ is equal to
the number of representatives of non-zero vectors in $N$ modulo
the action of $\{\pm 1\}$. This is also the number of
one-dimensional subspaces $F\E_v$ permuted by $\Sp_M$. This shows
that $\X\cong \ind_{\S_{v_n}}^{\Sp_M} F\E_{v_n}$, and completes
the proof of the theorem.\qed

If $l\neq 2$ we then have
$$
X=X_0\oplus X^+\oplus X^-.
$$
This is a decomposition of $X$ into non-isomorphic irreducible
$\Sp_M$-modules. If $l=2$ then $X^+=X^-$, which is isomorphic to
the $\Sp_M$-module $X/(X^+\oplus X_0)$ via the isomorphism
$E^+_v\mapsto E_v+(X^+\oplus X_0)$.

Given $\k\in F_q^*$, consider the non-trivial linear character
$\lambda[\k]$ of $F_q^+$ defined by
$$\lambda[\k](\a)=\lambda(\k\a),\quad \a\in F_q.$$
By varying $\k$ we obtain all non-trivial linear characters of
$F_q^+$.

\begin{thm}
\label{kk} Let $\k\in F_q^*$. Then
$$
\X_\lambda\cong \X_{\lambda[\k]}\text{ as }F\Sp_M\text{-modules
}\Leftrightarrow \k\text{ is a square}.
$$
\end{thm}

\noindent {\it Proof.} Suppose that $\X_\lambda$ and
$\X_{\lambda[\k]}$ are isomorphic as $F\Sp_M$-modules. The linear
characters entering the restriction of $\X_\lambda$ to $\Sp^M$
are of the form
$$
g\mapsto \lambda(\langle gv,v\rangle),
$$
where $0\neq v\in N$. One the other hand, the linear characters
entering the restriction of $\X_{\lambda[\k]}$ to $\Sp^M$ are of
the form
$$
g\mapsto \lambda(\k\langle gw,w\rangle).
$$
where $0\neq w\in N$. From the given hypothesis, it follows that
for some $v,w$ non-zero
$$\lambda(\langle gv,v\rangle)=\lambda(\k\langle gw,w\rangle)$$
for all $g\in\Sp^M$. Reasoning as in the proof of Lemma \ref{fp}
we see that
$$
\langle u,v\rangle^2=\k\langle u,w\rangle^2
$$
for all $u\in M$. Choosing $u$ so that $\langle u,v\rangle=1$ we
infer that $1=\k\langle u,w\rangle^2$, that is, $\k$ is a square.

Suppose conversely that $\k=\t^2$. Consider the $F$-isomorphism,
say $f$, from $X_\lambda$ to $X_{\lambda[\k]}$ defined by
$$
f(E_v)=E_{\t^{-1}v}.
$$
Then given any $s\in \Sp^M$ and any $a\in\Sp^M$ we have
$$
f(saE_v)=f(\lambda(\langle sav,av\rangle) E_{av})=
\lambda(\langle sav,av\rangle)E_{\t^{-1}av}.
$$
On the other hand
$$
\begin{aligned} saf(E_v) & =saE_{\t^{-1}v}=\lambda(\k \langle
sa\t^{-1}v,a\t^{-1}v\rangle)E_{a\t^{-1}v} =\lambda(\k
\t^{-1}\t^{-1}\langle sav,av\rangle)E_{\t^{-1}av}\\ & =
\lambda(\langle sav,av\rangle)E_{t^{-1}av}.
\end{aligned}
$$
It follows that $f$ is an isomorphism of $F\Sp^M$-modules.

Now $X^+\oplus X_0$ is the $1$-eigenspace of $-1\in Z(\Sp)$
acting on $X$,  and $X_0$ is the fixed points of $\Sp_M$. If
$l\neq 2$ then $X^-$ is the $-1$-eigenspace of $-1\in Z(\Sp)$
acting on $X$. Whether $l=2$ or not, we infer that when $f$
preserves $X^\pm$, whence $\X_\lambda\cong \X_{\lambda[\k]}$. \qed


\section{The Weil representation of $\Sp$}

We assume in this section that $p$ is odd, and rely on \cite{Sz}
as a general reference. Let $Y=Y_\lambda$ be an $F$-vector space
having a basis $(\epsilon_v)_{v\in N}$ indexed by $N$. Let $H_0$
the group whose underlying set is $F_q\times V$, with
multiplication
$$(\alpha,v)(\beta,w)=(\alpha+\beta+\langle v,w\rangle, v+w), \quad \alpha,\beta\in F_q,\;v,w\in
V.$$

The function $J:H_0\to\GL(Y)$ given by
$$
J(\alpha,u+v)\epsilon_w=\lambda(\alpha+\langle
u,v+2w\rangle)\epsilon_{w+v}, \quad u\in M,v,w\in N,\alpha\in F_q
$$
defines an absolutely irreducible $F$-representation of $H_0$. The
symplectic group acts on $H_0$ via the second coordinate:
${}^g(\alpha,v)=(\alpha,gv)$ for $g\in\Sp$, $\alpha\in F_q$,
$v\in V$. For $g\in\Sp$, the conjugate representation $J^g$ is
similar to $J$. Moreover, there is one, and only one (except when
$(n,q)=(1,3)$ when there are three), representation
$P:\Sp\to\GL(Y)$ satisfying:
$$
P(g)J(h)P(g)^{-1}=J({}^g h),\quad g\in\Sp, h\in H_0.
$$
We refer to $P$ as the Weyl representation of $\Sp$ over $F$ of
type $\lambda$ (when $(n,q)=(1,3)$ \emph{the} Weil representation
is the one given below). Consider the linear character
$\theta:\Sp_M\to F^*$ defined as follows for $s\in\Sp^M$ and
$a\in\Sp_{M,N}$:
$$
\theta(sa)=\theta(a)=\left(\frac{\mathrm{det}\,
a\vert_N}{q}\right)=\begin{cases} 1 & \text{ if } \mathrm{det}\,
a\vert_N\in {F_q^*}^2,\\
-1 & \text{ otherwise }.
\end{cases}
$$
The Weil representation can be defined on the basis vectors
$\epsilon_v$ of $Y$ by means of:
\begin{equation}
\label{larga}
P(sa)\epsilon_v=\theta(sa)^n\chi_{av}(s)\epsilon_{av},\quad
s\in\Sp^M, a\in\Sp_{M,N},
\end{equation}
\begin{equation}
\label{larga2} P(\rho_{v_n,-1})\epsilon_v=\left(\sum_{\a\in
F_q}\lambda(\a^2)\right)^{-1} \sum_{\a\in
F_q}\lambda(\a^2)\epsilon_{v+\a v_n}.
\end{equation}
We remark that \cite{Sz} defines the Weil representation over a
finite extension of $\Q$ containing a non-trivial $p$-th root of
unity. With slight modification of the arguments, one can see
that (\ref{larga})-(\ref{larga2}) hold when $\Q$ is replaced by
any field of characteristic different from $p$. Particular care
needs to be taken in the characteristic 2 case.

Having into account that $\Sp_M$ is a maximal subgroup of $\Sp$
and that $\rho_{v_n,-1}$ is not in $\Sp_M$, the above formulae
suffice to determine $P$. While explicit formulae exist for other
$P(g)$, $g\notin \Sp_M$, the above description of
$P(\rho_{v_n,-1})$ will be enough for our purposes.

For $v\in N\setminus\{0\}$, set
$\epsilon_v^\pm=\epsilon_v\pm\epsilon_{-v}$. Write
$Y_0=F\epsilon_0$, let $Y^-=Y^-_\lambda$ be the span of the
$\epsilon^-_v$, and let $Y^+=Y^+_\lambda$ be the span of the
$\epsilon^+_v$. Set $G(\lambda)=\sum_{\a\in F_q} \lambda(\a^2)$.

From (\ref{larga2}) we get
\begin{equation}
\label{j1}
\rho_{v_n,-1}\e_0=G(\lambda)^{-1}(\e_0+\sum_{\a\in{\cal{T}}_q}\lambda(\a^2)\e^+_{\a
v_n}),
\end{equation}
\begin{equation}
\label{j2} \rho_{v_n,-1}\e^+_{v}=G(\lambda)^{-1}( \sum_{\a\in
F_q}\lambda(\a^2)\e^+_{v+\a v_n}),\quad v\in N,v\notin F_q v_n,
\end{equation}
\begin{equation}
\label{j3} \rho_{v_n,-1}\e^+_{\beta
v}=G(\lambda)^{-1}(2\l(\beta^2)\e_0+ \sum_{\a\in {\cal {T}}_q
}(\lambda((\a+\beta)^2)+\lambda((\a-\beta)^2))\e^+_{\a
v_n}),\quad \beta\in {\cal{T}}_q.
\end{equation}
We are particularly interested in the case when $l=2$. In this
case $G(\l)=1$ because $\l(\a^2)+\l((-\a)^2)=2\l(\a^2)=0$ for all
$\a\in{\cal{T}}_q$, leaving only the term $\l(0)=1$. Thus, the
above equations simplify as follows when $l=2$:
\begin{equation}
\label{s1}
\rho_{v_n,-1}\e_0=\e_0+\sum_{\a\in{\cal{T}}_q}\lambda(\a^2)\e^+_{\a
v_n},
\end{equation}
\begin{equation}
\label{s2} \rho_{v_n,-1}\e^+_{v}=\sum_{\a\in
F_q}\lambda(\a^2)\e^+_{v+\a v_n},\quad v\in N,v\notin F_q v_n,
\end{equation}
\begin{equation}
\label{s3} \rho_{v_n,-1}\e^+_{\beta v_n}=\sum_{\a\in
{\cal{T}}_q}(\lambda((\a+\beta)^2)+\lambda((\a-\beta)^2))\e^+_{\a
v_n},\quad \beta\in {\cal{T}}_q.
\end{equation}
We also note that $\theta$ is trivial when $l=2$. In this case
(\ref{larga}) yields
\begin{equation}
\label{s4} sa\epsilon_v^+=\chi_{av}(s)\epsilon^+_{av},\quad
s\in\Sp^M, a\in\Sp_{M,N}.
\end{equation}

\begin{thm} If $l\neq 2$ then $Y^+\oplus Y_0$ and $Y^-$ are absolutely
irreducible $F\Sp$-modules. If $l=2$ then $Y^+=Y^-$ is an
absolutely irreducible $F\Sp$-module and $Y/Y^+\cong Y^+\oplus
Fy$, where $y$ is fixed by $\Sp$. In any case, $Y^\pm\cong
\theta^n X^\pm$ as $F\Sp_M$-modules. If $l=2$ then actually
$Y^+=Y^-\cong_{\Sp_M} X^-=X^+$.

\end{thm}

\noindent{\it Proof.} The last two assertions follow by comparing
(\ref{rep}) with (\ref{larga}), and having into account that
$\theta$ is trivial when $l=2$.

Suppose next that $l\neq 2$. Then $Y^-$ and $Y^+\oplus Y_0$ are
eigenspaces for $-1\in Z(\Sp)$ acting on $Y$, and are therefore
$\Sp$-stable. As the restriction of $Y^-$ to $\Sp_M$, namely
$\theta^n X^-$, is absolutely irreducible, so must be $Y^-$. Now
$Y^+\oplus Y_0$ is a decomposition into non-isomorphic irreducible
$\Sp_M$-modules. From (\ref{j1}) we see that $Y_0$ is not
$\Sp$-stable. By means of (\ref{j3}), and using $l\neq 2$, we
deduce that $Y^+$ is not $\Sp$-stable either. It follows that
$Y^+\oplus Y_0$ is an irreducible $\Sp$-module. As this argument
works for $F$ or any extension thereof, $Y^+\oplus Y_0$ is
absolutely irreducible.

Suppose next that $l=2$. Then $Y^+\oplus Y_0$ is still an
eigenspace for $-1\in Z(\Sp)$ and hence $\Sp$-stable, while
$Y^-=Y^+$. From (\ref{s2}) and (\ref{s3}) we see that now $Y^+$
is $\Sp$-stable. Since $Y^+\cong_{\Sp_M} X^+$, the $F\Sp$-module
$Y^+$ is absolutely irreducible. Consider the $F$-isomorphism
$Y/Y^+\to Y^+\oplus Fy$ given by
$$
\epsilon_0+Y^+\to y,\quad \epsilon_v+Y^+\to\epsilon_v^+,\,v\neq 0.
$$
Using (\ref{s1})-(\ref{s4}) we see that the action of $\Sp$ is
preserved. This completes the proof.\qed

We refer to $Y^-$ (resp. $Y^+\oplus Y_0$) as the Weil module of
$\Sp$ over $F$ of type $\lambda$ and degree $(q^n-1)/2$ (resp.
$(q^n+1)/2$). As our arguments are valid for any field of
characteristic $l\neq p$ containing a non-trivial $p$-th root of
unity, the above easily yields the following result.

\begin{thm} If $l\neq 2,p$ then the $l$-modular reduction of the
complex Weil modules of $\Sp$ of type $\lambda$ and degrees
$(q^n-1)/2$ and $(q^n+1)/2$ remain irreducible. The 2-modular
reduction of the complex Weil module of $\Sp$ of type $\lambda$
and degree $(q^n-1)/2$ remains irreducible, and is a constituent
of that of degree $(q^n+1)/2$. The other constituent is the
trivial module.
\end{thm}

In regards to the number of different types of Weil modules we
have the following.

\begin{thm} There are exactly two isomorphsim types of Weil modules of
$\Sp$ over $F$ of degree $(q^n-1)/2$ (resp. $(q^n+1)/2$).
\end{thm}

\noindent{\it Proof.} If $\k\in F_q^*$ is a square then the types
$\lambda$ and $\lambda[\k]$ are isomorphic. Indeed, the
isomorphism $f$ used in the proof of Theorem \ref{kk} is easily
seen to preserve the action of $\rho_{v_n,-1}$. Conversely, by
restricting to $\Sp_M$ and applying Theorem \ref{kk}, we see that
if the types $\lambda$ and $\lambda[\k]$ are isomorphic then
$\k\in F_q^*$ must be a square.\qed

For future reference we record the following results.

\begin{lem} Suppose $l=2$. Then
\label{lum1} $$ c_n \e^+_{v_n}=\sum_{\a\in {\cal{T}}_q}
(\lambda(-2\a)+\lambda(2\a))\e^+_{\a v_n}.
$$
\end{lem}

\noindent{\it Proof.} Observe the identity
\begin{equation}
\label{put1} c_n=\rho_{-1,u_n}\rho_{-1,v_n}\rho_{-1,u_n}.
\end{equation}
As $\rho_{-1,u_n}\in \Sp^M$, formula (\ref{s4}) yields
\begin{equation}
\label{put2} \rho_{-1,u_n} \e^+_{v}=\lambda(-\a_n^2)\e^+_{v},\quad
v=\a_1y_1+\cdots +\a_n v_n.
\end{equation}
By (\ref{s3}) we also have
$$
\rho_{v_n,-1}\e^+_{v_n}=\sum_{\a\in
{\cal{T}}_q}(\lambda((\a+1)^2)+\lambda((\a-1)^2))\e^+_{\a v_n}.
$$
Therefore
$$
\begin{aligned}
c_n\e^+_{v_n}
&=\rho_{-1,u_n}\rho_{-1,v_n}\rho_{-1,u_n}\e^+_{v_n}\\
&=\rho_{-1,u_n}\rho_{-1,v_n}\l(-1)\e^+_{v_n}\\
&=\l(-1)\rho_{-1,u_n}\sum_{\a\in
{\cal{T}}_q}(\lambda((\a+1)^2)+\lambda((\a-1)^2))\e^+_{\a v_n}\\
&=\sum_{\a\in
{\cal{T}}_q}(\lambda((\a+1)^2)+\lambda((\a-1)^2))\lambda(-1-\a^2)\e^+_{\a
v_n}\\
&=\sum_{\a\in {\cal{T}}_q} (\lambda(2\a)+\lambda(-2\a))\e^+_{\a
v_n}.
\end{aligned}
$$
\qed

\begin{lem} Suppose $l=2$. Then
\label{lum2}
$$ c_{n-1}\e^+_{v_n}=\sum_{\a\in F_q} \e^+_{v_n+\a v_{n-1}}.
$$
\end{lem}

\noindent{\it Proof.} From (\ref{s2}) we have
$$
\rho_{-1,v_n}\e^+_{v_{n-1}}=\sum_{\a\in
F_q}\l(\a^2)\e^+_{v_{n-1}+\a v_n}.
$$
Making use of (\ref{put1}) and (\ref{put2}) we get
$$
\begin{aligned}
c_n\e^+_{v_{n-1}}&=\rho_{-1,u_n}\rho_{-1,v_n}\rho_{-1,u_n}\e^+_{v_{n-1}}\\
&=\rho_{-1,u_n}\rho_{-1,v_n}\e^+_{v_{n-1}}\\
&=\rho_{-1,u_n}\sum_{\a\in F_q}\l(\a^2)\e^+_{v_{n-1}+\a v_n}\\
&=\sum_{\a\in F_q}\l(\a^2)\l(-\a^2)\e^+_{v_{n-1}+\a v_n}\\
&=\sum_{\a\in F_q}\e^+_{v_{n-1}+\a v_n}.
\end{aligned}
$$
Since $c_{n-1}=w_{n-1}c_nw_{n-1}$, the above yields
$$
c_{n-1}\e^+_{v_n}=w_{n-1}c_nw_{n-1}\e^+_{v_{n}}=w_{n-1}c_n\e^+_{v_{n-1}}=w_{n-1}\sum_{\a\in
F_q}\e^+_{v_{n-1}+\a v_n}=\sum_{\a\in F_q}\e^+_{v_{n}+\a v_{n-1}}.
$$
\qed

\section{The Steinberg representation of $\Sp$ restricted to $\Sp_M$}

For an element $w$ of the Weyl group $W$ of $\Sp$, let $n_w$ be
an element of $W_2$ satisfying $n_wH=w$. Consider the element
$$
\bar{e}=\sum_{b\in B}b\sum_{w\in W} (-1)^{\ell(w)}\;n_w
$$
of the symplectic group algebra $F\Sp$. The right ideal
$\bar{I}=\bar{e}F\Sp$ is a right $F\Sp$-module, considered by
Steinberg in \cite{St}. So far we have dealt exclusively with
left modules. We wish to adhere to this convention with
Steinberg's representation as well. For this purpose, consider
the involution $x\mapsto \bar{x}$ of $F\Sp$ that fixes all
scalars and inverts all symplectic transformations. Set
$$
e=(\sum_{w\in W}(-1)^{\ell(w)}\; n_w)\sum_{b\in B}b,
$$
and let $I=F\Sp\cdot e$. Then $I$ is naturally a left
$F\Sp$-module. Note that $\bar{I}$ is also a left $F\Sp$-module
under the rule
$$
g\cdot (\bar{e}\bar{x})=(\bar{e}\bar{x})\bar{g},\quad x\in
F\Sp,g\in \Sp.
$$
Furthermore, $I$ and $\bar{I}$ are isomorphic as left modules, via
the isomorphism $xe=\bar{e}\bar{x}$, $x\in F\Sp$.

We shall henceforth work with $e$ and $I$, and refer to
the latter as the (left) Steinberg module for $\Sp$ over $F$.
Our general reference for this section is \cite{St}.

An $F$-basis for $I$ is given by the $|U|$-elements $ue$, $u\in
U$. Thus $U$ affords the regular representation of $U$. The
derived quotient $U/U'$ is isomorphic to the direct product of
the $n$ fundamental root subgroups, and it is therefore an
elementary abelian $p$-group of order $q^n$. It follows that $F$
affords all linear characters of $U$. We infer that, given a
linear character $\sigma$ of $U$, there is a unique -up to
multiplication by a non-zero scalar- element $e_\sigma$ in $I$
upon which $U$ acts via $\sigma$. We may take
$$
e_\sigma=\sum_{u\in U}\sigma(u)^{-1}ue.
$$
Then
\begin{equation}
\label{u}
ue_\sigma=\sigma(u)e_\sigma,\quad u\in U.
\end{equation}

Let $w(i)$ be the element of $\Sp$ that preserves
the given symplectic basis, and satisfies:
$$
w(i)=(u_n,u_i,u_{i+1},...,u_{n-1})(v_n,v_i,v_{i+1},...,v_{n-1}).
$$
For $1\leq i\leq n$ and $\a\in F_q^*$, write $I_{\a,i}$ for the
$FU$-module generated by $w(i)e_{\chi_{\a v_n}}$.

\begin{thm}
\label{impo1} For
 $1\leq i\leq n$ and $\a\in F_q^*$ we have
$$
I_{\a,i}\cong_{FU}X^+_{\a,i}.
$$
\end{thm}

\noindent{\it Proof.} Observe first that the element $w(i)$ just
defined is the inverse of the element $w$, as defined in Theorem
\ref{irru}.

Next note that for $u\in U_w^+$ we have
$$
uw(i)e_{\chi_{\a v_n}}=w(i)w(i)^{-1}uw(i)e_{\chi_{\a v_n}}=
w(i)wuw^{-1}e_{\chi_{\a v_n}}.
$$
If $u\in U_w^+$ then ${}^w u\in U$, so (\ref{u}) applies, yielding
\begin{equation}
\label{uo}
uw(i)e_{\chi_{\a v_n}}=\chi_{\a v_n}({}^w u)w(i)e_{\chi_{\a v_n}}=
\chi_{\a v_n}(u^{w(i)})w(i)e_{\chi_{\a v_n}}\in I_{\a,i}.
\end{equation}

Since $U=U^+_w U^-_w$, it follows that $I_{\a,i}$ is generated
by the elements $uw(i)e_{\chi_{\a v_n}}$, as $u$ runs through
$U_w^-$. We proceed to show that these elements are linearly
independent.

Let $s\in\Sp^M$. For $u\in U^-_w$ we have
$$
suw(i)e_{\chi_{\a v_n}}=uu^{-1}suw(i)e_{\chi_{\a v_n}},
$$
where $s^u\in\Sp^M\subset U^+_w$. Hence by (\ref{uo}) and Lemma
\ref{g0}
$$
suw(i)e_{\chi_{\a v_n}}=u\chi_{\a v_n}(s^{uw(i)})w(i)e_{\chi_{\a
v_n}}=\chi_{\a u v_i}(s)u w(i)e_{\chi_{\a v_n}}.
$$
Thus $\Sp^M$ acts on $uw(i)e_{\chi_{\a v_n}}$, $u\in U_w^-$, via
the linear character $\chi_{\a u v_i}$. We contend these linear
characters are all different. Indeed, as noted in the proof of
Theorem \ref{irru}, the group $U^-_w$ consists of all $u\in U$
whose matrix has the form (\ref{b}), where $A$ is a lower
triangular with 1's on the diagonal, and all columns of $A$
different from column $i$ have zero entries below the diagonal.
It follows from this matrix description that $u\in U^-_w$ is
completely determined by what is does to $v_i$. In particular, we
see that $u\neq v$ in $U_w^-$ implies $\a u v_i\neq \pm \a v
v_i$. We deduce from Lemma \ref{fp} that $\chi_{\a u v_i}$ and
$\chi_{\a v v_i}$ are different linear characters of $\Sp^M$ for
$u\neq v\in U_w^-$. All in all, we get that the $uw(i)e_{\chi_{\a
v_n}}$, $u\in U_w^-$, must be linearly independent. As the also
generate $I_{\a,i}$, the form a basis of $I_{\a,i}$.

The preceding discussion shows that $I_{\a,i}$ affords a monomial
representation of $U$; the stabilizer of $Fw(i)e_{\chi_{\a v_n}}$
is $U^+_w$; $I_{\a,i}$ is the direct sum of left translates of
$Fw(i)e_{\chi_{\a v_n}}$ by elements of $U^-_w$, which is a
transversal for $U^+_w$ in $U$. We conclude $I_{\a,i}\cong
\ind_{U^+_w}^U Fw(i)e_{\chi_{\a v_n}}$ as $FU$-modules. By
(\ref{uo}) and Lemma \ref{g0} $Fw(i)e_{\chi_{\a v_n}}$ affords
the linear character $\chi_{\a,v_i}$ of $U^+_w$. The result now
follows from Theorem \ref{irru}.\qed

We next observe
$$
he=e,\quad h\in H.
$$
As $H$ normalizes $U$, for each linear character $\sigma$ of $U$
we may consider the linear character ${}^h \sigma$ of $U$,
defined by ${}^h \sigma(u)=\sigma(u^h)$, $u\in U$. Thus, for all
linear characters $\sigma$ of $U$ and all $h\in H$
\begin{equation}
\label{hu} he_\sigma=h\sum_{u\in U}\sigma(u)^{-1}ue=\sum_{u\in U}
\sigma(u)^{-1} ({}^h u) he=\sum_{u\in U}\sigma(u^h)^{-1}ue=
e_{({}^h \sigma)}.
\end{equation}

For $1\leq i\leq n$, recall the element $w(i)$ defined prior to
Theorem \ref{impo1}. Let $h\in H$. As $W_2$ normalizes $H$ we have
$h^{w(i)}\in H$. Since $v_i$ us an eigenvector for $h$, we may
write $hv_i=h_iv_i$ for some $h_i\in F_q^*$. From Lemma \ref{g0}
we obtain the formula
\begin{equation}
\label{hwi}
{}^{h^{w(i)}} \chi_{\a v_n}=\chi_{\a h_i v_n}.
\end{equation}
From (\ref{hu}) and (\ref{hwi}) we see that for $1\leq i\leq n$,
$h\in H$, and $\a\in F_q^*$
\begin{equation}
\label{trian} hw(i)e_{\chi_{\a v_n}}=w(i)h^{w(i)}e_{\chi_{\a
v_n}}=w(i)e_{\chi_{\a h_i v_n}}.
\end{equation}
It follows that for $1\leq i\leq n$, $h\in H$
and $\a\in F_q^*$
\begin{equation}
\label{iH} h I_{\a,i}=I_{h_i \a,i}.
\end{equation}

Now $\chi_{\a v_n}=\chi_{-\a,v_n}$, whence $I_{\a,i}=I_{-\a,i}$
for all $\a\in F_q^*$. As the $FU$-modules $I_{\a,i}$, $\a\in
{\cal{T}}_q$, are irreducible and non-isomorphic (cf. Theorems
\ref{irru}, \ref{isu} and \ref{impo1}), they are in direct sum
within $I$. Set $I_i=\oplus_{\a\in{\cal{T}}_q} I_{\a,i}$. From
(\ref{iH}) we see that $I_i$ is an $FB$-module, which is clearly
isomorphic to $X_i^+$.

In order to prove the next result we shall use, for the first time
in the paper, a beautiful identity due Steinberg, namely formula
(16) of \cite{St}.

\begin{thm}
\label{mejor} Let $w\in\{w_1,...,w_{n-1}\}$ and $v\in N$. Suppose
$U_w^+\subset S_v$. Then
$$
w(\sum_{u\in U^+_w}\sum_{u'\in U^-_w}\chi_{v}(u)^{-1}uu'
e)=\sum_{u\in U^+_w}\sum_{u'\in U^-_w}\chi_{w(v)}(u)^{-1}uu' e-
(q+1)\sum_{u\in U_w^+}\chi_{w(v)}(u)^{-1}ue.
$$
\end{thm}

\noindent{\it Proof.} We have
$$
w(\sum_{u\in U^+_w}\sum_{u'\in U^-_w}\chi_{v}(u)^{-1}uu' e)=
\sum_{u\in U^+_w}\sum_{u'\in U^-_w}\chi(u)^{-1}({}^w u)wu' e.
$$
Since $w$ has order 2, it normalizes $U_w^+$. Thus, it follows
from $U^+_w\subset S_v$ that $U^+_w\subset S_{w(v)}$. All in all,
$u\mapsto \chi_v({}^w u)$ is a linear character of $U_w^+$, which
by Lemma \ref{g0} must be equal to $\chi_{w(v)}$. We may thus
write
$$
\begin{aligned}
w(\sum_{u\in U^+_w}\sum_{u'\in U^-_w}\chi_{v}(u)^{-1}uu' e) &=
\sum_{u\in U^+_w}\sum_{u'\in U^-_w}\chi_v(u)^{-1}({}^w u)wu' e\\
&=\sum_{u\in U^+_w}\sum_{u'\in U^-_w}\chi_v(u^w)^{-1}uwu' e\\
&=\sum_{u\in U^+_w}\sum_{u'\in U^-_w}\chi_{w(v)}(u)^{-1}uwu' e.
\end{aligned}
$$
Since $we=-e$, to $u'=1$ there corresponds the summand
$$
- \sum_{u\in U_w^+}\chi_{w(v)}(u)^{-1}ue.
$$
As $U^-_w=X_{(i,i+1)}$ for some $i\in\{1,...,n-1\}$, we may write
all $u'\in U_w^-$ in the form
$$
u'=x_{i,i+1}(\a),\quad \a\in F_q.
$$
By means of Steinberg's relation (16) of \cite{St} we may write
$$
wx_{i,i+1}(\a)e=x_{i,i+1}(-\a^{-1})e-e,\quad \a\in F_q^*.
$$
It follows that $w(\sum_{u\in U^+_w}\sum_{u'\in
U^-_w}\chi_{v}(u)^{-1}uu' e)$ is equal to
$$
- \sum_{u\in U_w^+}\chi_{w(v)}(u)^{-1}ue -(q-1)\sum_{u\in
U_w^+}\chi_{w(v)}(u)^{-1}ue+ \sum_{u\in U^+_w}\sum_{\a\in
F_q^*}\chi_{w(v)}(u)^{-1}ux_{i,i+1}(-\a^{-1})e.
$$
Since $-\a^{-1}$ runs through $F_q^*$ as $\a$ does, we may
replace $-\a^{-1}$ by $\a$ in the above expression. Thus, by
adding and subtracting $\sum_{u\in U_w^+}\chi_{w(v)}(u)^{-1}ue$
we obtain the desired result.\qed

\begin{thm}
\label{refo} For all $\a\in F_q$ we have
$$w_ie_{\chi_{\a v_n}}=e_{\chi_{\a v_n}}-
(q+1)\sum_{u\in U_w^+}\chi_{\a v_n}(u)^{-1}ue,\quad
i\in\{1,...,n-2\},
$$
whereas
$$w_{n-1}e_{\chi_{\a v_n}}=\sum_{u\in U^+_w}\sum_{u' \in U^-_w}\chi_{\a v_{n-1}}(u)^{-1}uu' e-
(q+1)\sum_{u\in U_w^+}\chi_{\a v_{n-1}}(u)^{-1}ue.
$$
\end{thm}

\noindent{\it Proof.} Let $w\in\{w_1,...,w_{n-1}\}$. We know
that, not only $U_w^+\subset S_{\a v_n}$, but in fact $U$ is
included in $S_{\a v}$ and $\chi_{\a v_n}$ is a linear character
of $U$. As such, its kernel contains $U_w^-=X_{(i,i+1)}$ due to
Corollary \ref{rooto}. Thus
$$
e_{\chi_{\a v_n}}=\sum_{u\in U^+_w}\sum_{u'\in
U^-_w}\chi(uu')^{-1}uu' e=\sum_{u\in U^+_w} \sum_{u'\in
U^-_w}\chi(u)^{-1}uu' e.
$$
Since $w_i(\a v_n)=\a v_n$ when $i\in\{1,...,n-2\}$ and $w_n(\a
v_n)=\a v_{n-1}$, the result follows from Theorem \ref{mejor}.\qed

\begin{thm}
\label{great}
 Let $g\in\Sp_M$ and let $\a\in F_q^*$.
Then $ge_{\chi_{\a v_n}}=e_{\chi_{\a v_n}}$ implies $gv_n=\pm
v_n$. The converse holds precisely when $l$ divides $q+1$ (or
$n=1$).
\end{thm}

\noindent{\it Proof.} Suppose $g\in\Sp_M$ satisfies $ge_{\chi_{\a
v_n}}=e_{\chi_{\a v_n}}$. Given $s\in\Sp^M$ Lemma \ref{g0} gives
$$
sge_{\chi_{\a v_n}}=gs^g e_{\chi_{\a v_n}}=\chi_{\a
v_n}(s^g)ge_{\chi_{\a v_n}}=\chi_{\a gv_n}(s) e_{\chi_{\a v_n}}.
$$
On the other hand $se_{\chi_{\a v_n}}=\chi_{\a v_n}(s)e_{\chi_{\a
v_n}}$. Thus the linear characters $\chi_{\a g v_n}$ and
$\chi_{\a v_n}$ of $\Sp^M$ are equal. It follows from Lemma
\ref{fp} that $gv_n=\pm v_n$.

Suppose conversely that $gv_n=\pm v_n$. We wish to analyze under
what conditions $ge_{\chi_{\a v_n}}=e_{\chi_{\a v_n}}$. For ease
of notation we shall write $\chi=\chi_{\a v_n}$.

As $-1\in Z(\Sp)$ acts trivially on $I$, we we may assume without
loss of generality that $g v_n=v_n$. Write $g=sa$, where
$s\in\Sp^M$ and $a\in\Sp_{M,N}$. Now
$$
v_n=gv_n=sav_n=av_n+(sav_n-av_n).
$$
Here $av_n\in N$ and $sav_n-av_n\in M$. As $M$ and $N$ intersect
trivially, we infer $av_n=v_n$ and $sv_n=v_n$. Since $\Sp^M$ acts
on $e_\chi$ via $\chi$ and $sv_n=v_n$, it follows from
(\ref{formu}) that $se_\chi=\chi$. We may thus assume $g=a$.

By means of the usual $BN$-pair decomposition of $\Sp_{M,N}\cong
\GL(M)$ we may write $a=bwu$, where $b\in B$, $w\in W_0$, and
$u\in T$. Now $uv_n=v_n$, and therefore $ue_\chi=e_\chi$, so we
may dispense with $u$. From $bw v_n=v_n$ we easily see that
$bv_n=v_n$ and $w v_n=v_n$. Write $b=hv$, where $h\in H$ and
$v\in U$. Then $vv_n=v_n$ and $ve_\chi=e_\chi$. Hence $hv_n=v_n$,
and therefore $he_\chi=e_\chi$ by (\ref{trian}). Thus we may
assume $g=w$. As $w v_n=v_n$, $w\in\langle
w_1,..,w_{n-2}\rangle$. We may further suppose that $w=w_i$,
$i=1,...,n-2$. By Theorem \ref{refo} we have
$$we_\chi=e_\chi-(q+1)\sum_{u\in U_w^+}\chi(u)^{-1}ue.$$
It follows that $we_\chi=e_\chi$ if and only if $l$ divides
$q+1$.\qed

\begin{cor}
\label{ele} Suppose $l$ divides $q+1$. Let $\a\in F_q^*$ and
$w_1,w_2\in W_0$. Then $w_1e_{\chi_{\a v_n}}=w_2e_{\chi_{\a
v_n}}$ if and only if $w_1v_n=w_2v_n$.
\end{cor}

Let $L=L_\l$ be the sum of all $I_{i}$ inside $I$. As the $I_i$
are not isomorphic $FB$-modules, the sum is direct. It follows
that $\mathrm{dim}\,L=\mathrm{X^+}=|\Sp_M:\widehat{S}_{v_n}|$.

\begin{thm}
\label{cuatro} Suppose $l$ divides $q+1$. Then $L$ is an
absolutely irreducible $\Sp_M$-module isomorphic to $X^+_\lambda$.
\end{thm}

\noindent{\it Proof.} We first show that $L$ is
$\Sp_{M,N}$-invariant, affording a monomial representation.

A typical basis element of $L$ is of the form $uw(i)e_{\chi_{\a
v_n}}$ where $1\leq i\leq n$ and $u\in U_{{w(i)}^{-1}}^-$.
Multiplying this on the left by $g\in \Sp_{M,N}$ yields
$guw(i)e_{\chi_{\a v_n}}$, where $guw(i)$ is in $\Sp_{M,N}$, and
hence is of the form $bwv$. Here $b\in B$, $w\in W_0$ and $v\in
T$. We know that $v\in T$ acts trivially on $e_{\chi_{\a v_n}}$.
By Corollary \ref{ele} $we_{\chi_{\a v_n}}=w(j)e_{\chi_{\a
v_n}}$, where $wv_n=v_j$. Write $b=uh$, where $h\in H$ and $u\in
U$. Then $hw(j)e_{\chi_{\a v_n}}=w(j)e_{\chi_{\a h_i v_n}}$ by
(\ref{trian}), where $hv_j=h_j v_j$. Writing $u=u_1u_2$, where
$u_1\in U^-_{w(j)^{-1}}$ and $u_2\in U^+_{w(j)^{-1}}$, we see that
$uw(j)e_{\chi_{\a h_i v_n}}$, and hence $guw(i)e_{\chi_{\a
v_n}}$, is a scalar multiple of $u_1w(j)e_{\chi_{\a h_i v_n}}$.
Since this another typical basis element, we have shown that $L$
is $\Sp_{M,N}$-invariant and that $L$ affords a monomial
representation of $\Sp_{M,N}$ relative to the given basis
vectors, as claimed.

As seen in the proof of Theorem \ref{impo1}, $\Sp^M$ sends every
typical basis element to a scalar multiple of itself. We deduce
that $L$ is $\Sp_M$-invariant, and permutes the one dimensional
subspaces $Fuw(i)e_{\chi_{\a v_n}}$; the above argument shows
that the stabilizer of $Fe_{\chi_{v_n}}$ is $\widehat{S}_{v_n}$;
we have $\mathrm{dim}\,L=|\Sp_M:\widehat{S}_{v_n}|$. Therefore
$L\cong \ind_{\widehat{S}_{v_n}}^{\Sp_M} F e_{\chi_{v_n}}$. Since
the linear character of  $\widehat{S}_{v_n}$ afforded by
$Fe_{\chi_{v_n}}$ is $\chi^+_{v_n}$, the result follows from
Theorem \ref{irru}.\qed

Suppose $l$ divides $q+1$. If in the definition (\ref{formu}) of
$\chi_v$ we use $\l[\k]$, $\k\in F_q^*$, instead of $\l$, then $L$
will contain a copy of $X^+_{\l[\k]}$. Now for $\k$ a square
$L_\lambda=L_{\l[\k]}$. However if $\k$ is not a square then
$L_\l\not\cong L_{\l[\k]}$ by Theorems \ref{kk} and \ref{cuatro}.
We may combine these comments and the preceding results in the
following theorem.

\begin{thm} Suppose $l$ divides $q+1$.
The Steinberg module $I$ contains an irreducible
$F\Sp_M$-submodule isomorphic to $X^+_\l$. If $p$ is odd $I$ also
contains a copy of $X^+_{\l[\k]}$, where $\k\in F_q^*$ is not a
square. They are not isomorphic, so $I$ contains their direct
sum. Moreover, $X^+_\l$ and $X^+_{\l[\k]}$ are isomorphic to
$\theta^n Y_\l^+$ and $\theta^n Y_{\l[\k]}^+$, where $Y_\l$ and
$Y_{\l[\k]}$ are the Weil modules of types $\l$ and $\l[k]$
restricted to $\Sp_M$. If $l=2$ then $X^+_\l$ and $X^+_{\l[\k]}$
are isomorphic to the restriction to $\Sp_M$ of the Weil modules
$Y_\l^-$ and $Y_{\l[\k]}-$ of types $\l$ and $\l[k]$  and degree
$(q^n-1)/2$.
\end{thm}


\section{The Steinberg representation of $\Sp$}

We assume henceforth that $l=2$ and $p$ is odd. For the remainder
of the paper we shall write $\chi=\chi_{v_n}$. Let
$e_1=\sum_{u\in U}ue$. Then $S=Fe_1$ is the socle of
$F\Sp$-module $I$ (cf. Theorem 4.7 of \cite{G}), affording the
trivial representation of $\Sp$. Set $\LH=\LH_\l=L_\l\oplus S$ and
$\Lt=\Lt_\l=\LH_\l/S$. Our main goal is to show that the
$F$-subspace $\Lt_\l$ of $I/S$ is $\Sp$-stable and isomorphic to
the Weil module $Y^+_\lambda$ of degree $(q^n-1)/2$ as an
$F\Sp$-module. Our main tool will be the following criterion.

\begin{thm}
\label{nukes} Suppose the following identities hold in $I$:
\begin{equation}
\label{toon1} c_n e_\chi=e_1+\sum_{\a\in {\cal{T}}_q}
(\lambda(-2\a)+\lambda(2\a)) e_{\chi_{\a v_n}},
\end{equation}
\begin{equation}
\label{toon2} c_{n-1} e_\chi=\sum_{\a\in F_q}
w_{n-1}x_{n-1,n}(-\a) w_{n-1} e_\chi.
\end{equation}
Then $\Lt$ is $\Sp$-stable and $\Lt\cong Y^+$ as $F\Sp$-modules.
\end{thm}

\noindent{\it Proof.}  As $L$ and $Y^+$ are absolutely irreducible
isomorphic $\Sp_M$-modules, there is a unique
$F\Sp_M$-isomorphism $f:L\to Y^+$ up to multiplication by a
non-zero scalar. The dimension of the $F$-subspace in $L$ and $Y$
where $U$ acts via $\chi_{v_n}$ is equal to one. Thus, we may
assume without loss of generality that
\begin{equation}
\label{uno23} f(e_{\chi_{v_n}})=\e^+_{v_n}.
\end{equation}
Let $\f:\Lt\to Y^+$ be the $F\Sp_M$-isomorphism inherited from
$f$. Thus
$$
\f(x+S)=f(x),\quad x\in L.
$$
We wish to show that

(a) $\Lt$ is $\Sp$-stable;

(b) $\f$ is an $F\Sp$-isomorphism.

\noindent We first turn our attention to (a). We know that $\Sp_M$
and $c_n$ generate $\Sp$, as $\Sp_M$ is a maximal subgroup of
$\Sp$ and $c_n\notin \Sp_M$. Hence, to see that $\Lt$ is
$\Sp$-stable, we only need to verify that $\Lt$ is stable under
$c_n$.

A typical $F$-generator of $L$ is of the form $uhwe_\chi$ where
$u\in T$, $w\in W_0$ and $h\in H$. We need to show that
$c_n(uhwe_\chi)\in\LH$. We divide the proof into two cases.

Case 1. $wv_n=v_n$. In this case $c_n$ and $w$ commute. Since $T$
is contained in $U_{c_n}^+$, $c_n$ conjugates $u$ into $U$. As $H$
is normalized by $c_n$, it follows that ${}^{c_n} (uhw)$ belongs
to $\Sp_M$. Therefore ${}^{c_n} (uhw)L=L$, and a fortiori
${}^{c_n} (uhw)\LH=\LH$. Since $c_ne_\chi$ belongs to $\LH$ by
(\ref{toon1}), we infer
$$
c_n(uhwe_\chi)={}^{c_n} (uhw)c_ne_\chi\in\LH.
$$

Case 2. $wv_n=v_i$, $i\neq n$. By Corollary \ref{ele} we may write
$we_\chi= w'w_{n-1}e_\chi$ where $w'$ is an element of $W_0$ that
fixes $v_n$ and moves $v_{n-1}$ to $v_i$. By means of the identity
$w_{n-1}c_nw_{n-1}=c_{n-1}$ we may write
$$
c_n (uhw e_\chi)={}^{c_n}(uhw')c_nw_{n-1}e_\chi=
{}^{c_n}(uhw')w_{n-1}w_{n-1}
c_nw_{n-1}e_\chi={}^{c_n}(uhw')w_{n-1}c_{n-1}e_\chi.
$$

Now ${}^{c_n}u\in U$, ${}^{c_n}w'=w'$, and ${}^{c_n } h\in H$.
Therefore ${}^{c_n}(uhw')w_{n-1}\in\Sp_M$, and hence
${}^{c_n}(uhw')w_{n-1}\LH=\LH$. Since $c_{n-1}e_\chi\in\LH$ by
(\ref{toon2}), it follows that $c_n (uhw e_\chi)\in\LH$. This
proves (a).

In order to establish (b), it suffices to see that $\f$ commutes
with the action of $c_n$ on $\LH$, as $\f$ is already an
isomorphism of $F\Sp_M$-modules.

From Lemmas \ref{lum1} and \ref{lum2} we have
\begin{equation}
\label{cnn} c_n \e^+_{v_n}=\sum_{\a\in {\cal{T}}_q}
(\lambda(-2\a)+\lambda(2\a))\e^+_{\a v_n},
\end{equation}
\begin{equation}
\label{cnn2} c_{n-1}\e^+_{v_n}=\sum_{\a\in F_q} \e^+_{v_n+\a
v_{n-1}}.
\end{equation}

Let $\a\in F_q^*$, and let $h\in H$ satisfy $hv_n=\a v_n$. Then
$he_{\chi_{v_n}}=e_{\chi_{\a v_n}}$ by (\ref{trian}). Therefore
\begin{equation}
\label{ch} f(e_{\chi_{\a v_n}})=f(h
e_{\chi_{v_n}})=hf(e_{\chi_{v_n}})=h\e^+_{v_n}=\e^+_{\a v_n},\quad
\a\in F_q^*.
\end{equation}
It follows that
$$
\begin{aligned}
\f(c_n e_\chi+S) &= f(\sum_{\a\in {\cal{T}}_q}
(\lambda(-2\a)+\lambda(2\a))e_{\chi_{\a v_n}})\quad \text{ by
}(\ref{toon1})\\
&=\sum_{\a\in {\cal{T}}_q}
(\lambda(-2\a)+\lambda(2\a))f(e_{\chi_{\a v_n}})\\
&=\sum_{\a\in {\cal{T}}_q} (\lambda(-2\a)+\lambda(2\a))\e^+_{\a
v_n}\quad \text{ by
}(\ref{ch})\\
&= c_n \e^+_{v_n}\quad \text{ by
}(\ref{cnn})\\
&= c_n\f(e_\chi+S).
\end{aligned}
$$
We also have
$$
\begin{aligned}
\f(c_{n-1} e_\chi+S) &= f(\sum_{\a\in F_q} w_{n-1}x_{n-1,n}(-\a)
w_{n-1}
e_{\chi_{v_n}}) \text{ by }(\ref{toon2})\\
&=\sum_{\a\in F_q} w_{n-1}x_{n-1,n}(-\a) w_{n-1} f(e_{\chi_{v_n}})\\
&=\sum_{\a\in F_q} w_{n-1}x_{n-1,n}(-\a) w_{n-1}\e^+_{v_n}\\
&=\sum_{\a\in F_q} w_{n-1}x_{n-1,n}(-\a) \e^+_{v_{n-1}} \text{ by }(\ref{s4})\\
&=\sum_{\a\in F_q} w_{n-1}\e^+_{v_{n-1}+\a v_n}\text{ by }(\ref{s4})\\
&=\sum_{\a\in F_q} \e^+_{v_{n}+\a v_{n-1}}\text{ by }(\ref{s4})\\
&=c_{n-1}\e^+_{v_n}\text{ by }(\ref{cnn2})\\
&=c_{n-1}\f(e_\chi+S).
\end{aligned}
$$
Thus $\f$ commutes with the actions of $c_n$ and $c_{n-1}$ on the
single element $e_\chi$. We proceed to verify that this suffices
for $\f$ to commute with the action of $c_n$ on the whole $\Lt$.

Let $uhwe_\chi$ be a typical $F$-generator of $L$, where $u\in
T$, $w\in W_0$ and $h\in H$. Suppose first $wv_n=v_n$. Then
$$
\begin{aligned}
\f(c_n(uhwe_\chi+S)) &= \f({}^{c_n}(uhw)c_ne_\chi+S))\\
&={}^{c_n}(uhw)\f(c_ne_\chi+S)\quad\text{ since }{}^{c_n}(uhw)\in \Sp_M\\
&={}^{c_n}(uhw)c_n\f(e_\chi+S)\quad\text{ as just proven}\\
&=c_nuhw\f(e_\chi+S)\\
&=c_n\f(uhwe_\chi+S)\quad\text{ since }uhw\in\Sp_M.
\end{aligned}
$$
If now $wv_n=v_i$ with $i\neq n$, we may write $w=w'w_{n-1}$ as
above. Then
$$
\begin{aligned}
\f(c_n(uhwe_\chi+S)) &=\f({}^{c_n}(uhw')w_{n-1}c_{n-1}e_\chi+S))\\
&={}^{c_n}(uhw')w_{n-1}\f(c_{n-1}e_\chi+S)\quad\text{ since }{}^{c_n}(uhw')w_{n-1}\in \Sp_M\\
&={}^{c_n}(uhw')w_{n-1}c_{n-1}\f(e_\chi+S)\quad\text{ as just proven}\\
&=c_n u h w'w_{n-1}\f(e_\chi+S)\\
&=c_nu h w\f(e_\chi+S)\\
&=c_n\f(uhwe_\chi+S).
\end{aligned}
$$
Thus $\f$ commutes with the action of $c_n$ on all of $\Lt$.\qed

We next show that formula (\ref{toon1}) is indeed valid.

\begin{thm}
\label{unohipo} The following identity holds in $I$:
\begin{equation}
\label{uija} c_n e_\chi=e_1+\sum_{\b\in {\cal{T}}_q}
(\lambda(-2\b)+\lambda(2\b)) e_{\chi_{\b v_n}}.
\end{equation}
\end{thm}

\noindent {\it Proof.} We start by proving the identity
\begin{equation}
\label{tri} c_n e_\chi= \sum_{u\in U_{c_n}^+}\sum_{\a\in F_q^*}
\l(\a^{-1})u x_{n,2n}(\a)e.
\end{equation}
We know that $\chi$ is trivial on $U_{c_n}^+$ (cf. Lemma
\ref{puo}) and that $U_{c_n}^-$ is the root subgroup $X_{n,2n}$.
Moreover, $c_ne=-e$. Thus
$$
\begin{aligned}
c_n e_\chi & =c_n(\sum_{u\in U_{c_n}^+}\sum_{v\in U_{c_n}^-}
\chi(v)^{-1}u ve)\\
&=\sum_{u\in U_{c_n}^+}\sum_{\a\in F_q^*}
\chi(x_{n,2n}(\a))^{-1}({}^{c_n}u) c_n x_{n,2n}(\a)e- \sum_{u\in
U_{c_n}^+}({}^{c_n}u) e.
\end{aligned}
$$
In our context, Steinberg's formula (16) of \cite{St} reads
$$
c_n x_{n,2n}(\a)e=x_{n,2n}(-\a^{-1})e-e,\quad \a\in F_q^*.
$$
In the proof of Lemma \ref{puo} we established
$$
\chi(x_{n,2n}(\a))=\lambda(\a),\quad \a\in F_q.
$$
Therefore $\sum_{u\in U_{c_n}^+}\sum_{\a\in
F_q^*}\chi(x_{n,2n}(\a))^{-1}({}^{c_n}u) c_n x_{n,2n}(\a)e$ is
equal to
$$
\sum_{u\in U_{c_n}^+}\sum_{\a\in F_q^*} \l(-\a)({}^{c_n}u)
x_{n,2n}(-\a^{-1})e-\sum_{u\in U_{c_n}^+}\sum_{\a\in F_q^*}
\l(-\a)({}^{c_n})ue.
$$
As $\l$ is a non-trivial character  $F_q^*$
$$
\sum_{\a\in F_q^*}\l(-\a)=-\l(0)=-1.
$$
Thus
$$
-\sum_{u\in U_{c_n}^+}\sum_{\a\in F_q^*} \l(-\a)({}^{c_n}u)e= -
(\sum_{\a\in F_q^*}\l(-\a))(\sum_{u\in U_{c_n}^+}({}^{c_n}u) e)=
\sum_{u\in U_{c_n}^+}({}^{c_n}u) e.
$$
Since $c_n^2\in H$, conjugation by $c_n$ is an automorphism of
$U_{c_n}^+$. Hence
$$
\sum_{u\in U_{c_n}^+}({}^{c_n}u) e=\sum_{u\in U_{c_n}^+}u e,
$$
and
$$
\begin{aligned}
\sum_{u\in U_{c_n}^+}\sum_{\a\in F_q^*} \l(-\a){}^{c_n}u
x_{n,2n}(-\a^{-1})e &=\sum_{u\in U_{c_n}^+}\sum_{\a\in F_q^*}
\l(-\a){}u x_{n,2n}(-\a^{-1})e\\
&=\sum_{u\in U_{c_n}^+}\sum_{\a\in F_q^*} \l(\a^{-1}){}u
x_{n,2n}(\a)e.
\end{aligned}
$$
Combining the preceding equations we obtain (\ref{tri}).

In order to prove that the right hand sides of (\ref{uija}) and
(\ref{tri}) are equal we shall compare their coordinates relative
to the basis $(ue)_{u\in U}$ of $I$. A typical basis element is
of the form $ux_{n,2n}(\a)e$, where $u\in U_{c_n}^+$ and $\a\in
F_q$. The coefficient of $ux_{n,2n}(\a)e$ in the right hand side
of (\ref{tri}) is equal to 0 if $\a=0$, and $\l(\a^{-1})$
otherwise. Now the coefficient of $ux_{n,2n}(\a)e$ in the right
hand side of (\ref{uija}), say $C$, is equal to
$$
1+\sum_{\b\in {\cal{T}}_q} (\lambda(-2\b)+\lambda(2\b))\chi_{\b
v_n}(ux_{n,2n}(\a))^{-1}.
$$
Here
$$
\chi_{\b v_n}(ux_{n,2n}(\a))^{-1}=\chi_{\b
v_n}(x_{n,2n}(\a))^{-1}=\l(-\a\b^2).
$$
Thus
$$
C= 1+\sum_{\b\in {\cal{T}}_q}
(\lambda(-2\b)+\lambda(2\b))\l(-\a\b^2).
$$
If $\a=0$ then $C=1+\sum_{\b\in F_q^*}\l(2\b)=1+(-1)=0$, since
$\l$ is non-trivial and $p$ is odd. Suppose next that $\a\neq 0$.
We may write $C$ in the form
$$
1+\sum_{\b\in F_q^*}\lambda(2\b-\a\b^2).
$$
Suppose $\b$ is different from $2\a^{-1}$ and $\a^{-1}$. Then the
element $\gamma=2\a^{-1}-\b$ of $F_q$ is different from $0$ and
$\b$, and satisfies $2\b-\a\b^2=2\gamma-\a\gamma^2$. As $F$ has
characteristic 2, the summands corresponding to $\b$ and $\gamma$
yield
$$
\l(2\b-\a\b^2)+\l(2\gamma-\a\gamma^2)=2\l(2\b-\a\b^2)=0.
$$
Now when $\b=2\a^{-1}$, we have
$\lambda(2\b-\a\b^2)=\lambda(0)=1$. Added to the 1 on last formula
for $C$ yields 0 in $F$. Thus the only contributing summand for
$C$ is the one corresponding to $\b=a^{-1}$. It gives
$\lambda(2\a^{-1}-\a(\a^{-1})^2) =\lambda(\a^{-1})$. This
completes the proof.\qed

We proceed to verify the second condition of Theorem \ref{nukes}.
We shall do this by comparing the coefficient of each  basis
element $ue$, $u\in U$, on both sides of (\ref{toon2}). We begin
by considering the right hand side.

\begin{thm} The following identity holds for all $\a$ in $F_q$:
\label{op}
$$
w_{n-1}x_{n-1,n}(-\a)w_{n-1}e_\chi=\sum_{u\in
U_{w_{n-1}}^+}\sum_{v\in U_{w_{n-1}}^-}\chi_{v_{n}+\a
v_{n-1}}(u)^{-1}uve.
$$
\end{thm}

\noindent {\it Proof.} By virtue of Theorem \ref{refo} we have
$$
w_{n-1}e_\chi=\sum_{u\in U_{w_{n-1}}^+}\sum_{v\in U_{w_{n-1}}^-}
\chi_{v_{n-1}}(u)^{-1}uve.
$$
Therefore
$$
x_{n-1,n}(-\a)w_{n-1}e_\chi =\sum_{u\in U_{w_{n-1}}^+}\sum_{v\in
U_{w_{n-1}}^-}
\chi_{v_{n-1}}(u)^{-1}({}^{x_{n-1,n}(-\a)}u)x_{n-1,n}(-\a)ve.
$$
Observe that $X_{(n-1,n)}$ normalizes $U_{w_{n-1}}^+$. Thus
$u\mapsto \chi_{v_{n-1}}(u^{x_{n-1,n}(-\a)})$ is a linear
character of $U_{w_{n-1}}^+$, which by Lemma \ref{g0} is equal to
$\chi_{v_{n-1}+\a v_n}$. This makes sense, since $U_{w_{n-1}}^+$
fixes $v_n$ and $v_{n-1}$ modulo $M$, and it is therefore
contained in $S_{v_{n-1}+\a v_n}$. Thus
$$
\begin{aligned}
x_{n-1,n}(-\a)w_{n-1}e_\chi & = \sum_{u\in
U_{w_{n-1}}^+}\sum_{v\in U_{w_{n-1}}^-}
\chi_{v_{n-1}}(u^{x_{n-1,n}(-\a)})^{-1}u(x_{n-1,n}(-\a)v)e\\
&= \sum_{u\in U_{w_{n-1}}^+}\sum_{v\in U_{w_{n-1}}^-}
\chi_{v_{n-1}+\a v_n}(u)^{-1}u(x_{n-1,n}(-\a)v)e.
\end{aligned}
$$
As $v$ runs through $U_{w_{n-1}}^-$ so does $x_{n-1,n}(-\a)v$.
Hence
$$
x_{n-1,n}(-\a)w_{n-1}e_\chi = \sum_{u\in U_{w_{n-1}}^+}\sum_{v\in
U_{w_{n-1}}^-} \chi_{v_{n-1}+\a v_n}(u)^{-1}uve.
$$
We may now apply Theorem \ref{mejor} to obtain the desired
result.\qed

We wish to transform the above identity into something that can
later be compared with the left hand side of (\ref{toon2}). For
simplicity of notation we introduce the symplectic transformations
$E(b,c,d)\in \Sp^M$, defined by
$$
E(b,c,d)=x_{n-1,2n-1}(b)x_{n-1,2n}(c)x_{n,2n}(d)
$$
for all $b,c,d\in F_q$. For $a\in F_q$, we shall also write
$$
D(a)=x_{n-1,n}(a).
$$

Observe that every element of $U_{w_{n-1}}^+$ can be uniquely
written in the form $u E(b,c,d)$, where $u\in U_{c_n c_{n-1}}^+$
and $b,c,d\in F_q$. From Lemma \ref{puo} we infer that for such an
element
$$
\chi_{v_n+\a v_{n-1}}(uE(b,c,d))=\chi_{v_n+\a v_{n-1}}(E(b,c,d)).
$$
It follows from Theorem \ref{op} that
\begin{equation}
\label{unamas}
w_{n-1}D(-\a)w_{n-1}e_\chi= \sum_{u\in U_{c_n
c_{n-1}}^+}\sum_{d\in F_q} \sum_{c\in F_q} \sum_{b\in F_q}
\sum_{a\in F_q} \chi_{v_n+\a v_{n-1}}(E(b,c,d))^{-1} u
E(b,c,d)D(a)e.
\end{equation}

\begin{thm}
\label{pius} Let $u\in U_{c_n c_{n-1}}^+$ and let $a,b,c,d\in
F_q$. Let $C$ be the coefficient of $uE(b,c,d)D(a)e$ in the right
hand side of (\ref{toon2}). Then
$$
C=\begin{cases}
0 & \text{ if } b=0\text{ and }c\neq 0,\\
\lambda(-d) & \text{ if }b=0\text{ and }c=0,\\
\lambda(\frac{c^2}{b}-d) & \text{ if } b\neq 0.
\end{cases}
$$
\end{thm}

\noindent {\it Proof.} By (\ref{unamas})
$$
C=\sum_{\a\in F_q} \chi_{v_n+\a v_{n-1}}(E(b,c,d))^{-1}.
$$
From (\ref{formu})
$$
\chi_{v_n+\a v_{n-1}}(E(b,c,d))=\lambda(\langle E(b,c,d)(v_n+\a
v_{n-1}),v_n+\a v_{n-1}\rangle).
$$
Now
$$
\begin{aligned}
\langle E(b,c,d)(v_n+\a v_{n-1}),v_n+\a v_{n-1}\rangle &= \langle
c u_{n-1} + d u_n + \a b u_{n-1} + \a c u_n, \a
v_{n-1}+v_n\rangle\\
&= b\a^2+2c\a +d.
\end{aligned}
$$
It follows that
$$
C=\sum_{\a\in F_q} \lambda(-(b\a^2+2c\a+d))=\lambda(-d)\sum_{\a\in
F_q} \lambda(-(b\a ^2+2c\a)).
$$
If $b=0$ and $c=0$ the $C$ is equal to $q\lambda(-d)$, which
equals $\lambda(-d)$ in $F$ (as $q$ is odd and $l=2$). If $b=0$
and $c\neq 0$ then $-2c\a$ runs through $F_q$ as $\a$ runs
through $F_q$. Since $\lambda$ is non-trivial, $C$ equals 0.
Suppose next that $b\neq 0$. We inquire for which $\a$ in $F_q$
there exists $\b$ in $F_q$, different from $\a$, so that
$b\a^2+2c\a=b\b^2+2c\b$. This occurs precisely when $\a\neq
-\frac{c}{b}$. Since $l=2$ it follows that the only surviving
summand in the above formula for $C$ is
$$
\lambda(-(b  (-\frac{c}{b})^2+
2c(-\frac{c}{b}))=\lambda(\frac{c^2}{b}).
$$
Multiplying this by $\lambda(-d)$ we obtain the desired
result.\qed

We now turn our attention to the left hand side of (\ref{toon2}).

\begin{thm} The following relation holds in $I$:
$$
c_{n-1}e_\chi=\sum_{u\in U_{c_{n-1}c_n}^+}\sum_{d\in F_q}
\sum_{c\in F_q}\sum_{b\in F_q}\sum_{a\in F_q} \l(-d) u E(0,0,d)
w_{n-1} c_{n} E(0,c,b) D(a)e.
$$
\end{thm}

\noindent {\it Proof.}  Since $U_{c_{n-1}}^-$ is generated by the
$D(a)$ and the $E(b,c,0)$, it follows from Lemma \ref{puo} that
$\chi$ is trivial on $U_{c_{n-1}}^-$. Therefore
$$
e_\chi=\sum_{u\in U_{c_{n-1}}^+} \sum_{v\in U_{c_{n-1}}^-}
\chi(u)^{-1} u ve.
$$
Hence
$$
c_{n-1} e_\chi=\sum_{u\in U_{c_{n-1}}^+} \sum_{v\in U_{c_{n-1}}^-}
\chi(u)^{-1}({}^{c_{n-1}} u)(c_{n-1} v)e.
$$
As $c_{n-1}$ has order 2 modulo $H$, it follows that conjugation
by $c_{n-1}$ is an automorphism of $U_{c_{n-1}}^+$. Therefore
$$
c_{n-1} e_\chi=\sum_{u\in U_{c_{n-1}}^+} \sum_{v\in U_{c_{n-1}}^-}
\chi(u^{c_{n-1}})^{-1} u(c_{n-1} v)e.
$$
But $c_{n-1} v_n=v_n$, so Lemma \ref{g0} gives
$$
c_{n-1} e_\chi=\sum_{u\in U_{c_{n-1}}^+} \sum_{v\in U_{c_{n-1}}^-}
\chi(u)^{-1} u(c_{n-1} v)e.
$$
Now every element of $U_{c_{n-1}}^+$ can be uniquely written in
the form $uE(0,0,d)$, where $u\in U_{c_{n-1}c_n}^+$ and $d\in
F_q$. Also, every element of $U_{c_{n-1}}^-$ can be uniquely
written in the form $E(b,c,0)D(a)$, where $a,b,c\in F_q$. Since
$\chi$ is trivial on $U_{c_{n-1}c_n}^+$, it follows that
\begin{equation}
\label{gol} c_{n-1} e_\chi=\sum_{u\in U_{c_{n-1}c_n}^+}\sum_{d\in
F_q} \sum_{c\in F_q}\sum_{b\in F_q}\sum_{a\in F_q} \l(-d) u
E(0,0,d) (c_{n-1} E(b,c,0) D(a))e.
\end{equation}
Recall that $c_{n-1}=w_{n-1}c_n w_{n-1}$. Hence
$$
c_{n-1} E(b,c,0)D(a)e=w_{n-1} c_n w_{n-1} E(b,c,0)D(a)e=w_{n-1}
c_{n} E(0,c,b) w_{n-1} D(a) e.
$$
If $a=0$ then $w_{n-1}D(a)e=-e$, whereas if $a\neq 0$ then
$$
w_{n-1} D(a) e=D(-a^{-1})e-e.
$$
Summing over all $a\in F_q$, the right hand side of (\ref{gol})
thus yields $1+(q-1)=q$ terms equal to
$$
-\sum_{u\in U_{c_{n-1}c_n}^+}\sum_{d\in F_q} \sum_{c\in
F_q}\sum_{b\in F_q} \l(-d) u E(0,0,d) w_{n-1} c_{n} E(0,c,b)e,
$$
and one term equal to
$$
\sum_{u\in U_{c_{n-1}c_n}^+}\sum_{d\in F_q} \sum_{c\in
F_q}\sum_{b\in F_q}\sum_{a\in F_q^*} \l(-d) u E(0,0,d) w_{n-1}
c_{n} E(0,c,b) D(-a^{-1})e.
$$
As $-q\equiv 1\mod l$ and $-a^{-1}$ runs through $F_q$ when $a$
runs through $F_q$, the result follows.\qed

In order to proceed we record a number of relations in $\Sp$.
$$
D(a) E(b,c,d) D(a)^{-1}=E(b+2ac+a^2d,c+ad,d).
$$
In particular,
$$
D(a) E(b,0,0) D(a)^{-1}=E(b,0,0),
$$
$$
D(a) E(0,c,0) D(a)^{-1}=E(2ac,c,0),
$$
$$
D(a) E(0,0,d) D(a)^{-1}=E(a^2d,ad,d).
$$
We also have
$$
c_n D(a) c_n^{-1}=E(0,a,0),
$$
$$
c_n E(0,c,0) c_n^{-1}=D(-c),
$$
$$
c_n E(b,0,0) c_n^{-1}=E(b,0,0),
$$
$$
w_{n-1} E(b,c,d) w_{n-1}^{-1}=E(d,c,b).
$$
In regards to Steinberg's formula (16) of \cite{St}, we have the
following:
$$
c_n E(0,0,d)e=E(0,0,-d^{-1})e-e,\quad d\neq 0
$$
and
$$
w_{n-1} D(a)e=D(-a^{-1})e-e,\quad a\neq 0.
$$
In what follows we shall make implicit use of the above formulae.
Observe that
$$
\begin{aligned}
w_{n-1}c_n E(0,c,b)D(a)e &= w_{n-1}c_n E(0,c,0)D(a)D(-a)E(0,0,b)D(a)e\\
&=w_{n-1}c_n E(0,c,0)D(a) E(a^2 b,-ab,b)e\\
&=w_{n-1}D(-c)E(0,a,0)c_n E(a^2 b,-ab,b)e\\
&=w_{n-1} D(-c) E(0,a,0)E(a^2b,0,0)D(ab)c_n E(0,0,b)e.
\end{aligned}
$$
For $b\neq 0$, set
$$
f_1(a,b,c,d)=\l(-d)E(0,0,d)w_{n-1} D(-c)
E(0,a,0)E(a^2b,0,0)D(ab)E(0,0,-b^{-1})e,
$$
and under no restrictions on $b$, we write
$$
f_2(a,b,c,d)=\l(-d)E(0,0,d)w_{n-1} D(-c) E(0,a,0)E(a^2b,0,0)D(ab)
e.
$$
Then $c_{n-1}e_\chi$ is equal to
$$
\sum_{u\in U_{c_{n-1}c_n}^+}\sum_{d\in F_q} \sum_{c\in
F_q}\sum_{b^*\in F_q}\sum_{a\in F_q}uf_1(a,b,c,d)+ \sum_{u\in
U_{c_{n-1}c_n}^+}\sum_{d\in F_q} \sum_{c\in F_q}\sum_{b\in
F_q}\sum_{a\in F_q}uf_2(a,b,c,d).
$$
Now $w_{n-1} D(-c) E(0,a,0)E(a^2b,0,0)D(ab)E(0,0,-b^{-1})e$ is
equal to
$$
\begin{aligned}
& E(0,a,a^2b-2ac)w_{n-1}D(ab-c)E(0,0,-b^{-1})e =\\
& E(0,a,a^2b-2ac)w_{n-1} E(-(ab-c)^2b^{-1},-(ab-c)b^{-1},-b^{-1})D(ab-c)e =\\
& E(0,a,a^2b-2ac)E(-b^{-1},-(ab-c)b^{-1},-(ab-c)^2b^{-1})w_{n-1}D(ab-c)e =\\
& E(-b^{-1},cb^{-1},-c^2b^{-1})w_{n-1} D(ab-c)e\\
\end{aligned}
$$
while $w_{n-1} D(-c) E(0,a,0)E(a^2b,0,0)D(ab) e$ is equal to
$$
w_{n-1}E(a^2b-2ac,a,0)D(ab-c)e=E(0,a,a^2b-2ac)w_{n-1}D(ab-c)e.
$$
Set
$$
g_1(a,b,c,d)=\l(-d)E(-b^{-1},cb^{-1},-c^2b^{-1}+d)D((c-ab)^{-1}),\quad
b\neq 0,c\neq ab,
$$
$$
g_2(a,b,c,d)=\l(-d)E(-b^{-1},cb^{-1},-c^2b^{-1}+d),\quad b\neq 0,
$$
$$
g_3(a,b,c,d)=\l(-d)E(0,a,a^2b-2ac+d)D((c-ab)^{-1}),\quad c\neq ab,
$$
$$
g_4(a,b,c,d)=\l(-d)E(0,a,a^2b-2ac+d).
$$
Then $c_{n-1}e_\chi$ is equal to
$$
\begin{aligned}
&\sum_{u\in U_{c_{n-1}c_n}^+}\sum_{d\in F_q} \sum_{c\neq
ab}\sum_{b\neq 0}\sum_{a\in F_q}ug_1(a,b,c,d)+ \sum_{u\in
U_{c_{n-1}c_n}^+}\sum_{d\in F_q} \sum_{c\in
F_q}\sum_{b\neq 0}\sum_{a\in F_q}ug_2(a,b,c,d)+\\
& \sum_{u\in U_{c_{n-1}c_n}^+}\sum_{d\in F_q} \sum_{c\neq
ab}\sum_{b\in F_q}\sum_{a\in F_q}ug_3(a,b,c,d)+ \sum_{u\in
U_{c_{n-1}c_n}^+}\sum_{d\in F_q} \sum_{c\in F_q}\sum_{b\in
F_q}\sum_{a\in F_q}ug_4(a,b,c,d).
\end{aligned}
$$
Let $u\in U_{c_{n-1}c_n}^+$, $v\in U_{w_{n-1}}^-$ and $a,b,c,d\in
F_q$. Then

\noindent (1) If $b\neq 0$ and $v\neq 1$ then the coefficient of
$u E(-b^{-1},cb^{-1},-c^2b^{-1}+d)v e$ in $c_{n-1}e_\chi$ is
equal to $\lambda(-d)$. Throughout this process, the triple
$(b,c,d)$ is transformed into the triple
$(-b^{-1},cb^{-1},d-c^2b^{-1})$. This transformation is
invertible, and $(-b^{-1},-cb^{-1},d-c^2b^{-1})$ is the triple
transformed into $(b,c,d)$. It follows that, whenever $b\neq 0$,
the coefficient of $u E(b,c,d)v e$ in $c_{n-1}e_\chi$ is equal to
$\lambda(c^2b^{-1}-d)$.

\noindent(2) If $b\neq 0$ then the coefficient of $u
E(-b^{-1},cb^{-1},-c^2b^{-1}+d)e$ in $c_{n-1}e_\chi$ is equal to
$q\lambda(-d)=\l(-d)$. Reasoning as above, we see that if $b\neq
0$ the coefficient of $u E(b,c,d)e$ in $c_{n-1}e_\chi$ is equal
to $\lambda(c^2b^{-1}-d)$.

\noindent(3) If $v\neq 1$ then the coefficient of $uE(0,0,d)v$ in
$c_{n-1}e_\chi$ is equal to $q\l(-d)=\l(-d)$.

\noindent(4) The coefficient of $uE(0,0,d)$ in $c_{n-1}e_\chi$ is
equal to $q^2\l(-d)=\l(-d)$.

\noindent(5) If $a\neq 0$, $v\neq 1$ and $x\in F_q$ then the
coefficient of $uE(0,a,x)v$ in $c_{n-1}e_\chi$ is equal to
$\sum_{d\in F_q}(-d)=0$.

\noindent(6) If $a\neq 0$ and $x\in F_q$ the coefficient of
$uE(0,a,x)$ in $c_{n-1}e_\chi$ is equal to $q\sum_{d\in
F_q}(-d)=0$.

By comparing the above with Theorem \ref{pius} we obtain

\begin{thm}
\label{ufff} The following identity holds in $I$:
$$
c_{n-1} e_\chi=\sum_{\a\in F_q} w_{n-1}x_{n-1,n}(-\a) w_{n-1}
e_\chi.
$$
\end{thm}

Taking into account Theorems \ref{nukes}, \ref{unohipo} and
\ref{ufff} our main result is proven. If we then replace $\l$ by
$\l[\k]$, $\k$ not a square, in (\ref{formu}), we may state our
result as follows.

\begin{thm}
\label{crack} Let $q$ be a power of an odd prime $p$. Let $I$ be
the Steinberg module for $\Sp_{2n}(q)$ over a field $F$ of
characteristic 2 containing a primitive $p$-th root of unity.
Write $S$ for the socle of $I$, affording the trivial
representation of $\Sp_{2n}(q)$. Then $I/S$ contains as
irreducible modules the two non-isomorphic Weil modules for
$\Sp_{2n}(q)$ over $F$ of degree $(q^n-1)/2$.
\end{thm}

We wish to end the paper by locating the two composition factors
just found, in terms of the filtration for the Steinberg module
over $F$ introduced by Gow in \cite{G}. Adopting his notation, we
have

\begin{thm} Suppose $q\equiv 1\mod 4$. Let $\k$ be the highest
power of 2 dividing $|\Sp_{2n}(q):B|$. Then
$\overline{I(\k)}/\overline{I(\k-1)}$ contains as irreducible
modules the 2-modular reductions of the two non-isomorphic
complex Weil modules for $\Sp_{2n}(q)$ of degree $(q^n-1)/2$.
\end{thm}

\noindent{\it Proof.} It is a matter of computing the 2-valuation
of $|\Sp_{2n}(q):P_J|$, where $P_J$ is the parabolic subgroup of
$\Sp_{2n}(q)$ associated to the subset $J$ of
$\{w_1,...,w_{n-1},c_n\}$ corresponding to $\chi_{v_n}$. By
Corollary \ref{rooto}, $J=\{c_n\}$. Therefore $P_J=B\cup Bc_nB$,
a disjoint union. Since every element of $Bc_nB$ can be written
uniquely in the form $bc_nv$, where $b\in B$ and $v\in U_{c_n}^-$,
it follows that $|P_J|=(q+1)|B|$. Consequently, the 2-valuation of
$|\Sp_{2n}(q):P_J|$ is equal to $\k$ minus the 2-valuation of
$q+1$. As $q\equiv 1\mod 4$, the highest power of 2 dividing $q+1$
is 1, as required.\qed

\end{document}